\documentclass[10pt]{article}
\usepackage{amsthm, amsmath, amsfonts,	xcolor,	hyperref,	enumerate,	geometry, multirow,	caption, orcidlink}
\numberwithin{equation}{section}
\newtheorem{theorem}{Theorem}[section]
\newtheorem{example}[theorem]{Example}
\newtheorem{corollary}[theorem]{Corollary}

\newtheorem{remark}[theorem]{Remark}

\title{\textbf{General Chen's first inequality and applications for Riemannian maps}}

\author{Ravindra Singh\orcidlink{0009-0009-1270-3831}\thanks{Department of Mathematics, Banaras Hindu University, Varanasi, Uttar Pradesh-221005, India.\\ E-mail: khandelrs@bhu.ac.in}, Kiran Meena\orcidlink{0000-0002-6959-5853}\thanks{Department of Mathematics, Indian Institute of Technology Jodhpur, Rajasthan-342030, India.\\ E-mail: kiranmeena@iitj.ac.in; kirankapishmeena@gmail.com}, Kapish Chand Meena\orcidlink{0000-0003-0182-8822}\thanks{Scientific Analysis Group, Defence Research and Development Organisation, Delhi-110054, India. \& Defence Laboratory, Defence Research and Development Organisation, Jodhpur, Rajasthan-342011, India.\\ E-mail: meenakapishchand@gmail.com}}

\date{}

\begin{document}
	\maketitle
	\begin{abstract}
		\noindent In this paper, we propose \textit{general Chen's first inequality} for Riemannian maps between Riemannian manifolds and demonstrate its equality and sharpness via non-trivial examples. We also utilize this general inequality by establishing Chen's first inequalities when the target spaces are generalized complex and generalized Sasakian space forms, including real, complex, real K\"ahler, Sasakian, Kenmotsu, cosymplectic, and almost $C(\alpha)$ space forms. In addition, we estimate $\delta$-invariants under various hypotheses on these space forms. Finally, we validate our new approach by comparing particular results with those of existing approaches.
	\end{abstract}
	
	\noindent \textbf{Keywords: }{Riemannian maps; Chen's first inequality; $\delta$-invariants; Complex \& contact structures; Generalized space forms}\\
	
	\noindent \textbf{MSC Classification:} {53B35, 53C15, 53D15}

	\section{Introduction}
	
	In differential geometry, there is a shortage of relations between the extrinsic and intrinsic properties of submanifolds. In 1993, Chen \cite{Chen_1993} introduced a potential tool for Riemannian submanifolds relating external and internal invariants, now known as \textit{Chen's first inequality}. Notably, it helps to bound and relate curvatures, establish rigidity theorems, explore geometric flows, etc. \cite{Chen_2011}. The concept of a Riemannian map also has rich geometric and physical implications, allowing us to compare the geometric properties of the source and target manifolds \cite{Fischer_1992, Sahin_book}. In \cite{Sahin_2016, Meena}, \c{S}ahin et al. investigated Chen's first inequality for Riemannian maps from Riemannian manifolds to real and complex space forms. They observed that it relates the curvatures of the horizontal and range spaces. This raises a natural question: how does Chen's first inequality extend to Riemannian maps into other space forms, as in classical submanifold theory (see \cite{ACKY, AEMMO, Chen_1993, Chen_1996, Chen_2023, GAH, Kim_Choe, Kim_Kim, KST, Ozgur_2011, Ozgur_De, MMT_2003, Yoon_2004})? Here, we address this question by introducing an effective solution, namely, \textit{general Chen's first inequality}. Specifically, it can be used to establish Chen's first inequality for Riemannian maps with other target space forms. For instance, we use it for generalized complex and generalized Sasakian space forms, including specific cases. 
    
    Throughout the manuscript, the notations are interpreted as follows:
    \begin{table}[h]
    \begin{tabular}{ll}
			\hline
            \textbf{Notation} &\textbf{Interpretation}\\
            \hline
			$R^M$ (resp. $R^N$)  &Riemannian curvature tensor of $M$ (resp. $N$) \\
			$\nabla^M$ (resp. $\nabla^N$) & Levi-Civita connection on $M$ (resp. $N$)\\
			$2\rho^{\mathcal{H}}$ (resp. $2\rho^{\mathcal{R}}$) & Scalar curvature of $\mathcal{H}$ (resp. $\mathcal{R}$)\\
			$K^{\mathcal{H}}(\mathbb{P})$ (resp. $K^{\mathcal{R}}(\mathbb{P})$) & Sectional curvature of plane spanned by $\{h_1, h_2\}$ (resp. $\{\pi_\ast h_1, \pi_\ast h_2\}$)\\
			$\tau^{\mathcal{H}}$ & Tension field of $\mathcal{H}$\\
			$\beta_{ij}^{k}$  &$g_2\left( \left( \nabla \pi_{\ast }\right) \left(h_{i},h_{j}\right), V_{k}\right)$\\
			$\|\mathcal{P}\|^{2}$ & $\sum \limits_{i,j=1}^{r}\left(g_{2}(\pi_{\ast }h_{i}, \mathcal{P} \pi_{\ast }h_{j})\right)^{2}$\\
			$\Theta(\mathbb{P})$ & $\left(g_{2}( \pi_{\ast }h_{1}, \mathcal{P} \pi_{\ast }h_{2})\right)^{2}$\\
			$\Phi(\mathbb{P})$ & $\left(\eta(\pi_\ast h_1)\right)^2 + \left(\eta(\pi_\ast h_2)\right)^2$\\
			\hline
		\end{tabular}
        \end{table}
	
	Let us state the aforementioned general inequality:
	\begin{theorem}[\textbf{General Chen's First Inequality}]\label{main_lemma}
		Let $\pi:\left(M^{m}, g_1\right) \to \left( N^{n},g_2\right)$ be a Riemannian map between Riemannian manifolds with rank $3 \leq r < n$. Then for all $p\in M$ and every plane section $\mathbb{P} \subset T_{p}M$ spanned by $\left\{h_{1},h_{2}\right\}$, we have
		\begin{equation}\label{eq-GSSF-(1)}
			K^{\mathcal{H}} (\mathbb{P}) \geq \frac{1}{2}\left\{2\rho^{\mathcal{H}} -\frac{(r-2)}{(r-1)} \|\tau^{\mathcal{H}}\|^2 - 2\rho^{\mathcal{R}} + 2 K^{\mathcal{R}} ( \mathbb{P} ) \right\},
		\end{equation}
		where $2 \rho^{\mathcal{H}}$ (resp. $2 \rho^{\mathcal{R}}$) is the scalar curvature of $\mathcal{H}$ (resp. $\mathcal{R}$). The equality in $(\ref{eq-GSSF-(1)})$ holds if and only if there exist orthonormal bases $\left\{h_{i}\right\}_{i=1}^{r}$ and $\left\{V_{k}\right\} _{k=r+1}^{n}$ of $\mathcal{H}$ and $\mathcal{R}^{\perp}$, respectively, such that the shape operator becomes
		\begin{align}\label{eq-GSSF-(2)}
			&\mathcal{S}_{r+1}=\left(
			\begin{array}{ccccc}
				\beta_{11}^{r+1} & 0 & 0 & \cdots & 0 \\ 
				0 & \beta_{22}^{r+1} & 0 & \cdots & 0 \\ 
				0 & 0 & \beta_{11}^{r+1}+\beta_{22}^{r+1} & \cdots & 0 \\ 
				\vdots & \vdots & \vdots & \ddots & \vdots \\ 
				0 & 0 & 0 & \cdots & \beta_{11}^{r+1}+\beta_{22}^{r+1}
			\end{array}
			\right)_{r \times r} \nonumber \\&~{\rm and}~\mathcal{S}_{\ell}=\left( 
			\begin{array}{ccccc}
				\beta_{11}^{\ell } & \beta_{12}^{\ell } & 0 & \cdots & 0 \\ 
				\beta_{12}^{\ell } & -\beta_{11}^{\ell } & 0 & \cdots & 0 \\ 
				0 & 0 & 0 & \cdots & 0 \\ 
				\vdots & \vdots & \vdots & \ddots & \vdots \\ 
				0 & 0 & 0 & \cdots & 0\\
			\end{array}
			\right)_{r \times r} {\rm for}~ r+2 \leq \ell \leq n.
		\end{align}
        Geometrically, in the equality case when $r < n-1$, the leaves of $\mathcal{R}$ are
        \begin{enumerate}[$(i)$]
            \item $2$-quasi-umbilical submanifolds of the target space, if $\beta_{11}^{r+1} \neq 0, \beta_{22}^{r+1} \neq 0, \beta_{11}^{r+1} \neq \beta_{22}^{r+1}$, and either $\beta_{11}^\ell \neq 0$ or $\beta_{12}^\ell \neq 0$ for all $r+2 \leq \ell \leq n$. 
            
            \item totally $2$-quasi-umbilical submanifolds of the target space, if $\beta_{11}^{r+1} = - \beta_{22}^{r+1} \neq 0$, and either $\beta_{11}^{\ell} = \pm \beta_{11}^{r+1} = \mp \beta_{22}^{r+1}$, $\beta_{12}^\ell =0$, or $\beta_{12}^{\ell} = \pm \beta_{11}^{r+1} = \mp \beta_{22}^{r+1}$, $\beta_{11}^\ell =0$ for all $r+2 \leq \ell \leq n$.

            \item totally geodesic submanifolds of the target space, if $\beta_{11}^{r+1} = \beta_{22}^{r+1} = 0$ and $\beta_{11}^\ell = \beta_{12}^\ell = 0$ for all $r+2 \leq \ell \leq n$.
        \end{enumerate}
        In addition, geometrically, in the equality case when $r = n-1$, the leaves of $\mathcal{R}$ are
        \begin{enumerate}[$(i)$]
            \item totally $2$-quasi-umbilical submanifolds of the target space, if $\beta_{11}^{r+1} \neq 0, \beta_{22}^{r+1} \neq 0$, and $\beta_{11}^{r+1} \neq \beta_{22}^{r+1}$.

            \item totally quasi-umbilical submanifolds of the target space, if exactly one condition in $\{\beta_{11}^{r+1} = 0, \beta_{22}^{r+1} = 0\}$ holds.

            \item totally geodesic submanifolds of the target space, if $\beta_{11}^{r+1} = \beta_{22}^{r+1} = 0$.
        \end{enumerate}
	\end{theorem}
	Clearly, general Chen's first inequality also relates the sectional, scalar, and squared mean curvatures of the horizontal and range spaces. By utilizing this general inequality, we establish the following Chen's first inequalities for Riemannian maps to generalized complex and generalized Sasakian space forms, respectively:
	
\begin{theorem}\label{Theorem_GCSF_1} 
		Let $\pi:\left( M^{m},g_{1}\right) \to \left(N^{n}\left( f_{1}, f_{2}\right), g_{2}, J\right)$ be a Riemannian map from a Riemannian manifold to a generalized complex space form with rank $3 \leq r < n$. Then for all $p\in M$ and every plane section $\mathbb{P} \subset T_{p}M$ spanned by $\left\{h_{1},h_{2}\right\}$, the following hold.
        \begin{enumerate}[$(1)$]
            \item If $f_1$ and $f_2$ are constant, then
            \begin{equation}\label{eq-GCSF-(1)}
			K^{\mathcal{H}}\left( \mathbb{P} \right) \geq \frac{1}{2}\left\{ 2\rho^{\mathcal{H}}-\frac{\left(r-2\right) }{\left( r-1\right) }\|\tau^{\mathcal{H}}\|^{2}- (r-2)(r+1) f_1 -3f_{2}\left( \|\mathcal{P}\|^{2}-2\Theta(\mathbb{P}) \right) \right\}.
		\end{equation}
       
        \item If $f_1$ and $f_2$ are non-constant, then as $n=4$ and $r=3$, we have
        \begin{equation}\label{const_case_gcsf}
			K^{\mathcal{H}}\left( \mathbb{P} \right) \geq \frac{1}{2}\left\{ 2\rho^{\mathcal{H}}-\frac{1}{2}\|\tau^{\mathcal{H}}\|^{2}- 4 f_1 -3f_{2}\left( \|\mathcal{P}\|^{2}-2\Theta(\mathbb{P}) \right) \right\}.
		\end{equation}
        \end{enumerate}
        In addition, the equality conditions are identical to those of $(\ref{eq-GSSF-(1)})$.
	\end{theorem}
	
	\begin{theorem}\label{theorem_GSSF}
		Let $\pi:\left( M^{m},g_{1}\right) \to \left( N^{n}\left(f_1, f_2, f_3\right), g_{2}, \phi, \xi, \eta\right)$ be a Riemannian map from a Riemannian manifold to a generalized Sasakian space form with rank $3 \leq r < n$. Then for all $p\in M$ and every plane section $\mathbb{P} \subset T_{p}M$ spanned by $\left\{h_{1},h_{2}\right\}$, the following hold.
		\begin{enumerate}[$(1)$]
			\item If Reeb vector $\xi \in \Gamma(\mathcal{R})$, then
			\begin{align}\label{eq-GNSSF-(1)}
				K^{\mathcal{H}}\left( \mathbb{P} \right) \geq & \frac{1}{2}\left\{ 2\rho^{\mathcal{H}}-\frac{\left(r-2\right) }{\left( r-1\right) }\|\tau^{\mathcal{H}}\|^{2}-f_{1}\left( r^{2}-r-2\right) \right. \nonumber\\& \left.-3f_{2}\left( \|\mathcal{P}\|^{2}-2~ \Theta(\mathbb{P})\right) + 2 f_3 (r-1-\Phi(\mathbb{P}))\right\}.
			\end{align}
			
			\item If Reeb vector $\xi \in \Gamma(\mathcal{R}^\perp)$, then
			\begin{equation}\label{CFI_when_xi_in_R_perp}
				K^{\mathcal{H}}\left( \mathbb{P} \right) \geq \frac{1}{2}\left\{ 2\rho^{\mathcal{H}}-\frac{\left(r-2\right) }{\left( r-1\right) }\|\tau^{\mathcal{H}}\|^{2}- (r-2)(r+1) f_1 -3f_{2}\left( \|\mathcal{P}\|^{2}-2 ~\Theta(\mathbb{P}) \right) \right\}.
			\end{equation}
		\end{enumerate}
		In addition, the equality conditions are identical to those of $(\ref{eq-GSSF-(1)})$.
	\end{theorem}
	
	This paper is assembled in eight sections. In Section \ref{sec_prelims}, we collect the preliminary information required for the study. Section \ref{sec_main_lemma} proves Theorem \ref{main_lemma}. In Section \ref{section_examples}, we provide two non-trivial examples that satisfy the inequality $(\ref{eq-GSSF-(1)})$. Further, Section \ref{sec_utilization} is dedicated to the proofs of Theorems \ref{Theorem_GCSF_1} and \ref{theorem_GSSF}. Section \ref{sec_corollary} covers Chen's first inequalities when the target spaces are real, complex, real K\"ahler, Sasakian, Kenmotsu, cosymplectic, and almost $C(\alpha)$-space forms. In Section \ref{sec_deltainv}, we estimate the $\delta$-invariants. The last section is dedicated to a concluding remark.
	
    \section{Prerequisites}\label{sec_prelims}
	
	In this section, we present some basic concepts to prove the main inequalities.
	
	Let $\pi: (M^m, g_1) \to (N^n, g_2)$ be a smooth map between Riemannian manifolds with $0 < {\textrm{\normalfont rank }} \pi < \min\{m,n\}$, and let $\pi_{\ast p}: T_p M \to T_{\pi(p)} N$ be its differential map at $p$. Denoting the kernel space of $\pi_\ast$ at $p \in M$ by $\mathcal{V}_p = (\ker \pi_{\ast p})$ and its orthogonal complementary space in the tangent space $T_p M$ by $\mathcal{H}_p = (\ker \pi_{\ast p})^\perp$, we have
	\begin{equation*}
		T_p M = \mathcal{V}_p \oplus \mathcal{H}_p.
	\end{equation*}
	Similarly, we have
	\begin{equation*}
		T_{\pi(p)} N  = \mathcal{R}_p \oplus \mathcal{R}_p^\perp,
	\end{equation*}
	where $\mathcal{R}_p = ({\textrm{\normalfont range }} \pi_{\ast p})$ denotes the range of $\pi_\ast$, and its orthogonal complementary space is $\mathcal{R}_p^\perp= ({\textrm{\normalfont range }} \pi_{\ast p})^\perp$ in the tangent space $T_{\pi(p)} N$. Then, the map $\pi$ is said to be a \textit{Riemannian map} if, for all $X, Y \in \Gamma(\mathcal{H})$, we have \cite{Fischer_1992}
	\begin{equation}\label{riemannian_map}
		g_1(X, Y) = g_2(\pi_\ast X, \pi_\ast Y).
	\end{equation}
	The map $\pi_\ast$ can be viewed as a section of the bundle $\operatorname{Hom}(TM,\pi^{-1}TN)$ $\to M$, where $\pi^{-1}TN$ is the pullback bundle. The bundle $\operatorname{Hom}(TM,\pi^{-1}TN)$ has a connection $\nabla$ induced by $\nabla^M$. With this, the symmetric \textit{second fundamental form} of $\pi$ is given by \cite{Nore_1986}
	\begin{equation}\label{eqn_sff}
		(\nabla \pi_\ast) (\vartheta_1,\vartheta_2) = \nabla_{\pi_\ast \vartheta_1}^N \pi_\ast \vartheta_2 - \pi_\ast({\nabla}_{{\vartheta_1}}^M \vartheta_2),
	\end{equation}
	for all $\vartheta_1, \vartheta_2 \in \Gamma(TM)$. In addition, for $X, Y \in \Gamma(\mathcal{H})$, we have $(\nabla \pi_\ast) (X,Y)\in \Gamma(\mathcal{R}^\perp)$ \cite{Sahin_2010}. Moreover, its trace defines the \textit{tension field} $\tau^{\mathcal{H}}$ of $\mathcal{H}$. In other words, suppose $\{h_i\}_{i=1}^{r}$ is an orthonormal basis of $\mathcal{H}$, then by \cite{Sahin_book}
	\begin{equation*}
		\tau^{\mathcal{H}} = \sum\limits_{i=1}^{r} ((\nabla \pi_{\ast})(h_i, h_i)).
	\end{equation*}
	
	For any vector field $X$ on $M$ and any section $V$ of $\mathcal{R}^\perp$, we have $\nabla_X^{\pi\bot} V$, which is the orthogonal projection of $\nabla_X^N V$ on $\mathcal{R}^\perp$, where $\nabla^{\pi \bot}$ is a linear connection on $\mathcal{R}^\perp$ such that $\nabla^{\pi \bot} g_2 = 0$. Then for $\pi$, the \textit{shape operator} $\mathcal{S}_V$ is defined as \cite[p. 188]{Sahin_book}: 
	\begin{equation*}
		\nabla_{\pi_\ast X}^N V = -\mathcal{S}_V \pi_\ast X + \nabla_X^{\pi \bot} V.
	\end{equation*}
	At $p\in M$, we have $\nabla_{\pi_\ast X}^N V(p) \in T_{\pi(p)} N$, $\mathcal{S}_V \pi_\ast X \in \pi_{\ast p} (T_p M)$, and $\nabla_{X}^{\pi \bot} V(p) \in (\pi_{\ast p}(T_p M))^\bot$.

    \begin{remark}\label{rem_on_rank}
		Throughout the paper, we assume a Riemannian map of rank $3 \leq r < n$ to avoid weak relations.
	\end{remark}
	
	An almost Hermitian manifold $(N, g_2, J)$ is said to be a generalized complex space form (denoted $N(f_1, f_2)$) if there exist smooth functions $f_1$ and $f_2$ on $N$ such that \cite{Olszak, Tricerri}
	\begin{align}\label{curvature_tensor_gcsf}
		R^N(\zeta_1, \zeta_2)\zeta_3 &= f_1\{g_2(\zeta_2, \zeta_3) \zeta_1 - g_2(\zeta_1, \zeta_3)\zeta_2\} \nonumber \\&+ f_2 \{g_2(\zeta_1, J \zeta_3)J\zeta_2 - g_2(\zeta_2, J\zeta_3) J\zeta_1 + 2g_2(\zeta_1, J\zeta_2)J\zeta_3\},
	\end{align}
	for all $\zeta_1, \zeta_2, \zeta_3 \in \Gamma(TN)$. In addition, an almost contact manifold $(N, g_2, \phi, \xi, \eta)$ is said to be a generalized Sasakian space form (denoted $N(f_1, f_2, f_3)$) if there exist smooth functions $f_1$, $f_2$, and $f_3$ on $N$ such that \cite{ABC}
	\begin{align}\label{curvature_tensor_gssf}
		&R^N(\zeta_1, \zeta_2)\zeta_3 = f_1 \{g_2(\zeta_2, \zeta_3)\zeta_1 - g_2(\zeta_1, \zeta_3)\zeta_2\}\nonumber\\&+f_2\{g_2(\zeta_1, \phi \zeta_3) \phi \zeta_2 - g_2(\zeta_2, \phi \zeta_3)\phi \zeta_1 + 2g_2(\zeta_1, \phi \zeta_2)\phi \zeta_3\}\nonumber\\&+f_3\{\eta(\zeta_1) \eta(\zeta_3) \zeta_2 - \eta(\zeta_2) \eta(\zeta_3) \zeta_1 + g_2(\zeta_1, \zeta_3) \eta(\zeta_2) \xi - g_2(\zeta_2,\zeta_3) \eta(\zeta_1)\xi\},
	\end{align}
	for all $\zeta_1, \zeta_2, \zeta_3 \in \Gamma(TN)$. \\
	
	For the sake of standardizing and minimizing notations in the proofs, for any vector $\zeta \in \Gamma(TN)$, we correspondingly take
    \begin{equation}\label{decomposition_with_J}
		J \zeta =  \mathcal{P} \zeta + \mathcal{Q} \zeta \quad \text{or} \quad \phi \zeta =  \mathcal{P} \zeta + \mathcal{Q} \zeta,
	\end{equation}
	where $\mathcal{P} \zeta \in \Gamma(\mathcal{R})$, $\mathcal{Q} \zeta \in \Gamma(\mathcal{R}^\perp)$, and $\mathcal{P}$ is the endomorphism on $\mathcal{R}$.
    
	\section{Proof of general Chen's first inequality}\label{sec_main_lemma}
	
	In this section, we present the proof of Theorem \ref{main_lemma} as follows:
	
	\begin{proof}[Proof of Theorem \ref{main_lemma}]
		For $X, Y, Z, W \in \Gamma(\mathcal{H})$, we have the Gauss equation \cite[p. 189]{Sahin_book}
        \begin{eqnarray*}
			g_1(R^{M}(X,Y)Z, W) &=& g_2(R^{N}(\pi_{\ast }X,\pi_{\ast}Y)\pi_{\ast }Z,\pi_{\ast }W)  \\&&+g_2(\left( \nabla \pi_{\ast }\right) (Y,Z),\left( \nabla \pi_{\ast}\right) (X,W))\\&&-g_2(\left( \nabla \pi_{\ast }\right)(X, Z),\left( \nabla \pi_{\ast }\right) (Y, W)).
		\end{eqnarray*}
        Let $\left\{ h_{i}\right\}_{i=1}^{r}$ and $\left\{ \pi_{\ast }h_{i}\right\}_{i=1}^{r}$ be orthonormal bases of $\mathcal{H}$ and $\mathcal{R}$, respectively. Then by contraction, we get
		\begin{eqnarray*}
            \sum \limits_{i,j=1}^{r} g_1(R^{M}(h_{i},h_{j})h_{j},h_{i}) &=& \sum \limits_{i,j=1}^{r} g_2(R^{N}(\pi_{\ast }h_{i},\pi_{\ast}h_{j})\pi_{\ast }h_{j},\pi_{\ast }h_{i})  \\&&+\sum \limits_{i,j=1}^{r}g_2(\left( \nabla \pi_{\ast }\right) (h_{i},h_{i}),\left( \nabla \pi_{\ast}\right) (h_{j},h_{j}))\\&&-\sum \limits_{i,j=1}^{r}g_2(\left( \nabla \pi_{\ast }\right)(h_{i},h_{j}),\left( \nabla \pi_{\ast }\right) (h_{j},h_{i})).
		\end{eqnarray*}
		This gives
		\begin{equation}\label{eq-GSSF-(5)}
		2\rho^{\mathcal{H}}=2\rho^{\mathcal{R}}+\|\tau^{\mathcal{H}}\|^{2} -\|\nabla \pi_{\ast }\|^{2},
		\end{equation}
		where $2\rho^\mathcal{H}$ (resp. $2\rho^{\mathcal{R}}$) denotes the \textit{scalar curvature} of $\mathcal{H}$ (resp. $\mathcal{R}$), and given by
		\begin{equation*}
			2 \rho^\mathcal{R} = \sum \limits_{i,j=1}^{r} g_2 (R^N(\pi_\ast h_i, \pi_\ast h_j) \pi_\ast h_j, \pi_\ast h_i) ~ \text{and}~2 \rho^\mathcal{H} = \sum \limits_{i,j=1}^{r} g_1 (R^M(h_i, h_j) h_j, h_i).
		\end{equation*}
		Suppose
		\begin{equation}\label{eq-GSSF-(6)}
			\varepsilon =2\rho^{\mathcal{H}}-2\rho^{\mathcal{R}}-\frac{\left( r-2\right) }{\left(r-1\right) }\|\tau^{\mathcal{H}}\|^{2}.
		\end{equation}
		Then (\ref{eq-GSSF-(5)}) and (\ref{eq-GSSF-(6)}) yield,
		\begin{equation*}
			\|\tau^{\mathcal{H}}\|^{2}=\left( r-1\right) \left( \varepsilon+\|\nabla \pi_{\ast }\|^{2}\right) .  
		\end{equation*}
		Let $\mathbb{P} = \operatorname{span}\left\{ h_{1},h_{2}\right\} $, and let $\left\{ V_{k}\right\} _{k=r+1}^{n}$ be an orthonormal basis of $\mathcal{R}^{\perp}$ such that $V_{r+1}$ is parallel to $\tau^{\mathcal{H}}$ and $\beta_{ij}^{\alpha }=$ $g_2\left( \left( \nabla \pi_{\ast }\right) \left(h_{i},h_{j}\right) ,V_{\alpha }\right)$. Then, we have
        \begin{equation}\label{extra_eq}
            \left( \sum_{i=1}^{r} \beta_{ii}^{r+1}\right)^{2}=\left( r-1\right) \left\{ \sum_{i=1}^{r}\left( \beta_{ii}^{r+1}\right)^{2}+\sum_{i\neq j=1}^{r}\left(\beta_{ij}^{r+1}\right)^{2}+\sum_{i,j=1}^{r}\sum_{\ell =r+2}^{n}\left(\beta_{ij}^{\ell }\right)^{2}+\varepsilon \right\}.
        \end{equation}
		Employing \cite[Lemma 3.1]{Chen_1993}, we get
		\[
		2\beta_{11}^{r+1}\beta_{22}^{r+1}\geq \sum_{i\neq j=1}^{r}\left( \beta_{ij}^{r+1}\right)^{2}+\sum_{i,j=1}^{r}\sum_{\ell =r+2}^{n}\left( \beta_{ij}^{\ell }\right)^{2}+\varepsilon. 
		\]
		 By performing some computations, it can be expressed as
		\begin{eqnarray}\label{eq-GSSF-(8)}
			\beta_{11}^{r+1}\beta_{22}^{r+1} &\geq &\left( \beta_{12}^{r+1}\right)^{2}+\sum_{j>2}^{r}\left( \left( \beta_{1j}^{r+1}\right)^{2}+\left(\beta_{2j}^{r+1}\right)^{2}\right) +\frac{1}{2}\sum_{i\neq j>2}^{r}\left(\left( \beta_{ij}^{r+1}\right)^{2}\right)  \nonumber \\&&+\sum_{\ell =r+2}^{n}\left( \beta_{12}^{\ell }\right)^{2}+\sum_{\ell=r+2}^{n}\sum_{j>2}^{r}\left( \left( \beta_{1j}^{\ell }\right)^{2}+\left(\beta_{2j}^{\ell }\right)^{2}\right) +\sum_{\ell=r+2}^{n}\sum_{i,j>2}^{r}\left( \beta_{ij}^{\ell }\right)^{2}  \nonumber \\&&+\sum_{\ell =r+2}^{n}\left( \frac{1}{2}\left( \beta_{11}^{\ell }+\beta_{22}^{\ell}\right)^{2}-\beta_{11}^{\ell }\beta_{22}^{\ell }\right) +\frac{\varepsilon }{2}.
		\end{eqnarray}
		Particularly, if $h_i=h_{1}$ and $h_j =h_{2}$, then by Gauss equation we have
		\[
		\beta_{11}^{r+1}\beta_{22}^{r+1}=K^{\mathcal{H}}(\mathbb{P}) -K^{\mathcal{R}}\left( \mathbb{P} \right)
		+\left( \beta_{12}^{r+1}\right)^{2}-\sum_{\ell =r+2}^{n}\left( \beta_{11}^{\ell
		}\beta_{22}^{\ell }-\left( \beta_{12}^{\ell }\right)^{2}\right), 
		\]
		where the \textit{sectional curvatures} are defined by
		\begin{equation*}
			K^\mathcal{H}(\mathbb{P}) = g_1(R^M(h_1, h_2) h_2, h_1) ~ \text{and}~K^\mathcal{R}(\mathbb{P}) = g_2(R^N(\pi_\ast h_1, \pi_\ast h_2) \pi_\ast h_2, \pi_\ast h_1).
		\end{equation*}
		By using the aforementioned equation in (\ref{eq-GSSF-(8)}), we get
		\begin{eqnarray*}
			K^{\mathcal{H}}(\mathbb{P}) &\geq &K^{\mathcal{R}}\left( \mathbb{P} \right) +\sum_{\ell=r+1}^{n}\sum_{j>2}^{r}\left( \left( \beta_{1j}^{\ell }\right)^{2}+\left(\beta_{2j}^{\ell }\right)^{2}\right) \\&&+\frac{1}{2}\sum_{i\neq j>2}^{r}\left( \left( \beta_{ij}^{r+1}\right)^{2}\right) +\sum_{\ell=r+2}^{n}\sum_{i,j>2}^{r}\left( \beta_{ij}^{\ell}\right)^{2}+\sum_{\ell =r+2}^{n}\frac{1}{2}\left(\beta_{11}^{\ell}+\beta_{22}^{\ell }\right)^{2}+\frac{\varepsilon }{2}.
		\end{eqnarray*}
		And hence,
		\begin{equation}\label{eq-GSSF-(9)}
			K^{\mathcal{H}}(\mathbb{P}) \geq \frac{\varepsilon }{2} + K^{\mathcal{R}} \left( \mathbb{P} \right).
		\end{equation}
		From (\ref{eq-GSSF-(6)}) and (\ref{eq-GSSF-(9)}) we obtain (\ref{eq-GSSF-(1)}). Further, the equality in (\ref{eq-GSSF-(9)}) holds if and only if
        \begin{equation*}
            \sum_{\ell=r+1}^{n}\sum_{j>2}^{r}\left( \left( \beta_{1j}^{\ell }\right)^{2}+\left(\beta_{2j}^{\ell }\right)^{2}\right) + \frac{1}{2}\sum_{i\neq j>2}^{r}\left( \left( \beta_{ij}^{r+1}\right)^{2}\right) +\sum_{\ell=r+2}^{n}\sum_{i,j>2}^{r}\left( \beta_{ij}^{\ell}\right)^{2}+\sum_{\ell =r+2}^{n}\frac{1}{2}\left(\beta_{11}^{\ell}+\beta_{22}^{\ell }\right)^{2} = 0.
        \end{equation*}
        This means\\
		$\beta_{1j}^{\ell} =\beta_{2j}^{\ell}=0 \ {\rm for\ all}\ j>2, \ell \in \{r+1, r+2, \dots, n\}; \\
		\beta_{ij}^{r+1} =0 \ {\rm for\ } i\neq j>2; \\
		\beta_{ij}^{\ell } =0\ {\rm for}\ i,j \in \{3,4, \dots, r\},\ \ell \in \{r+2, \dots, n\}; \\
		\beta_{11}^{\ell }+\beta_{22}^{\ell } =0\ {\rm for\ all}\ \ell \in \{r+2, \dots, n\}$. \\
		
		\noindent We can choose an orthonormal basis $\left\{h_{1},h_{2}\right\} $ such that $\beta_{12}^{r+1}=\beta_{21}^{r+1}=0$. Moreover, since (\ref{extra_eq}) holds, using \cite[Lemma 3.1]{Chen_1993}, we also have
        \begin{equation*}
        \beta_{11}^{r+1}+\beta_{22}^{r+1}=\beta_{33}^{r+1}=\cdots =\beta_{rr}^{r+1}.
        \end{equation*}
        Thus, $\mathcal{S}$ becomes as desired.\\
        
        Further when $r<n-1$, (\ref{eq-GSSF-(2)}) can be viewed as follows:
        \begin{itemize}
            \item In the direction $V_{r+1}$, the shape operator is diagonal and has trace $(r-1)(\beta_{11}^{r+1} + \beta_{22}^{r+1})$ with the eigenvalues $\beta_{11}^{r+1}$ (with multiplicity one), $\beta_{22}^{r+1}$ (with multiplicity one), and $\beta_{11}^{r+1} + \beta_{22}^{r+1}$ (with multiplicity $(r-2)$).

            \item In the directions $\{V_\ell\}_{\ell=r+2}^n$, the shape operator is traceless with the eigenvalues $\pm \sqrt{(\beta_{11}^\ell)^2 + (\beta_{12}^\ell)^2}$ (with multiplicities one) and $0$ (with multiplicity $(r-2)$).
        \end{itemize}
        Geometrically, in the equality case, the leaves of $\mathcal{R}$ with dimension $r<n-1$ have three distinct principal curvatures in every direction $\{V_k\}_{k=r+1}^{n}$ if $\beta_{11}^{r+1} \neq 0, \beta_{22}^{r+1} \neq 0, \beta_{11}^{r+1} \neq \beta_{22}^{r+1}$, and if either $\beta_{11}^\ell \neq 0$ or $\beta_{12}^\ell \neq 0$ for all $r+2 \leq \ell \leq n$. Hence, they are $2$-quasi-umbilical submanifolds of the target space \cite{Meulemeester_Verstraelen}. In particular, if $\beta_{11}^{r+1} =- \beta_{22}^{r+1} \neq 0$, and if either $\beta_{11}^{\ell} = \pm \beta_{11}^{r+1} = \mp \beta_{22}^{r+1}, \beta_{12}^\ell =0$ or $\beta_{12}^{\ell} = \pm \beta_{11}^{r+1} = \mp \beta_{22}^{r+1}, \beta_{11}^\ell =0$ for all $r+2 \leq \ell \leq n$, then the leaves of $\mathcal{R}$ have three distinct principal curvatures $\pm \beta_{11}^{r+1} = \mp \beta_{22}^{r+1} \neq 0$ (with multiplicities one) and $0$ (with multiplicity $(r-2)$) in every direction $\{V_k\}_{k=r+1}^{n}$. Hence, in this particular case, they are totally $2$-quasi-umbilical submanifolds of the target space \cite{Meulemeester_Verstraelen}. Moreover, if $\beta_{11}^{r+1} = \beta_{22}^{r+1} = 0$ and if $\beta_{11}^\ell = \beta_{12}^\ell = 0$ for all $r+2 \leq \ell \leq n$, then the leaves of $\mathcal{R}$ have single principal curvature, namely $0$ (with multiplicity $r$) in every direction $\{V_k\}_{k=r+1}^{n}$; and hence they are totally geodesic \cite{Chen_totally_umbilical}.\\ 
        
        On the other hand, when $r= n-1$, we have the single direction $V_{r+1}$, and by (\ref{eq-GSSF-(2)}), the shape operator is diagonal and has trace $(r-1)(\beta_{11}^{r+1} + \beta_{22}^{r+1})$ with the eigenvalues $\beta_{11}^{r+1}$ (with multiplicity one), $\beta_{22}^{r+1}$ (with multiplicity one), and $\beta_{11}^{r+1} + \beta_{22}^{r+1}$ (with multiplicity $(r-2)$). Further, similarly, one can prove the equality case by taking into account the required notions from \cite{Chen_totally_umbilical, Chen_Verstraelen, Meulemeester_Verstraelen}.
	\end{proof}

    \section{Applicability of general Chen's first inequality}\label{section_examples}
    
    In this section, we illustrate the strictness and equality of the general Chen's first inequality. Precisely, we provide non-trivial examples of Riemannian maps for which the general inequality holds either strictly or with equality.\\

    First, we give an example of a Riemannian map for which the inequality (\ref{eq-GSSF-(1)}) holds with equality.

    \begin{example}
    Let $M= \left\{ \left(x_{1},x_{2},x_{3},x_{4},x_{5},x_{6}\right) \in {\mathbb{R}}^{6} : x_{i} \neq 0 \right\}$ and $N =\left\{ \left( y_{1},y_{2},y_{3},y_{4}\right) \in {\Bbb R}^{4} : y_2 y_3 y_4 \neq 0\right\}$ be two Riemannian manifolds with Riemannian metrics $g_1 = x_{6}^{2} \sum \limits_{i=1}^{5}dx_{i}^{2}+dx_{6}^{2}$ and $g_2 = dy_{1}^{2}+y_{4}^{2}dy_{2}^{2}+y_{4}^{2}dy_{3}^{2}+dy_{4}^{2}$, respectively. Define a map $ \pi: \left( M, g_{1}\right) \to \left( N,g_{2}\right) $ by
    \[
    \pi \left( x_{1},x_{2},x_{3},x_{4},x_{5},x_{6}\right) =\left(0,x_{4},x_{5},x_{6}\right).
    \]
    Then we have
    \begin{equation*}
        \mathcal{H} = \operatorname{span} \{h_1 = e_4, h_2 = e_5, h_3 = e_6\}, \quad \mathcal{V} = \operatorname{span} \{e_1, e_2, e_3\},
    \end{equation*}
    \begin{equation*}
        \mathcal{R} =\operatorname{span} \left\{ \pi_{\ast } h_{1}=e_{2}^{\ast }, \pi_{\ast } h_{2} =e_{3}^{\ast }, \pi_{\ast} h_{3} =e_{4}^{\ast }\right\}, \quad \mathcal{R}^{\perp} = \operatorname{span} \left\{V_4 =e_{1}^{\ast }\right\},
    \end{equation*}
    where $\left\{ e_{1}=\frac{1}{x_{6}}\frac{\partial }{\partial x_{1}},e_{2}=\frac{1}{x_{6}}\frac{\partial }{\partial x_{2}},e_{3}=\frac{1}{x_{6}}\frac{\partial }{\partial x_{3}},e_{4}=\frac{1}{x_{6}}\frac{\partial }{\partial x_{4}},e_{5}=\frac{1}{x_{6}}\frac{\partial }{\partial x_{5}},e_{6}=\frac{\partial }{\partial x_{6}}\right\}$\\ and $\left\{ e_{1}^{\ast }=\frac{\partial }{\partial y_{1}},e_{2}^{\ast }=\frac{1}{y_{4}}\frac{\partial }{\partial y_{2}},e_{3}^{\ast }=\frac{1}{y_{4}}\frac{\partial }{\partial y_{3}},e_{4}^{\ast }=\frac{\partial }{\partial y_{4}}\right\}$ are bases of $T_p M$ and $T_{\pi (p)} N$, respectively. We observe that $(\ref{riemannian_map})$ holds for all $h_{i}, h_{j} \in \Gamma (\mathcal{H})$. Thus, $\pi$ is a Riemannian map with $r=3$. Now, we compute the non-zero Christoffel symbols for $g_1$ as
    \begin{equation*}
        \Gamma_{x_i x_i}^{x_6} = -x_6 \quad \text{and} \quad \Gamma_{x_i x_6}^{x_i} = \frac{1}{x_6}\quad \text{for \quad $1 \leq i \leq 5$},
    \end{equation*}
    and for $g_2$ as
    \begin{eqnarray*}
    \Gamma_{y_2 y_2}^{y_4} =-y_{4}, \quad \Gamma_{y_3 y_3}^{y_4}=-y_{4}, \quad
    \Gamma_{y_2 y_4}^{y_2} =\Gamma_{y_4 y_2}^{y_2}=\frac{1}{y_{4}}, \quad 
    \Gamma_{y_3 y_4}^{y_3} =\Gamma_{y_4 y_3}^{y_3}=\frac{1}{y_{4}}.
\end{eqnarray*}
Then, using $(\ref{eqn_sff})$, we obtain
\[
\left( \nabla \pi_{\ast }\right) \left( h_{i}, h_{j}\right) =0,\ {\rm for\ all\ } h_{i}, h_{j} \in \Gamma(\mathcal{H}).
\]
Consequently, $\beta_{ij}^4 = g_2((\nabla \pi_\ast) (h_i, h_j), V_4) = 0$ for all $h_{i}, h_{j} \in \Gamma(\mathcal{H})$. Further, we compute
\begin{equation*}
    \|\tau^{\mathcal{H}} \|^2 = 0,
\end{equation*}\begin{equation*}
    K^{\mathcal{R}} (\mathbb{P}) = g_2(R^N(\pi_\ast h_1, \pi_\ast h_2) \pi_\ast h_2, \pi_\ast h_1) = g_2(R^N(e_2^*, e_3^*)e_3^*,e_2^*) = -\left(\frac{1}{y_4}\right)^2,
\end{equation*}
\begin{equation*}
    2 \rho^{\mathcal{R}} = \sum\limits_{i, j= 1}^{3} g_2(R^N(\pi_\ast h_i, \pi_\ast h_j) \pi_\ast h_j, \pi_\ast h_i) = -2\left(\frac{1}{y_4}\right)^2,
\end{equation*}
\begin{equation*}
    K^{\mathcal{H}} (\mathbb{P}) = g_1(R^M(h_1, h_2) h_2, h_1) = g_1(R^M(e_4, e_5)e_5,e_4) = -\left(\frac{1}{x_6}\right)^2,
\end{equation*}
and
\begin{equation*}
    2 \rho^{\mathcal{H}} = \sum\limits_{i, j= 1}^{3} g_1(R^M(h_i, h_j) h_j, h_i) = -2\left(\frac{1}{x_6}\right)^2.
\end{equation*}
Thus, the inequality in $(\ref{eq-GSSF-(1)})$ holds with equality as we have $y_4 = x_6$ from the defined map $\pi$. This is also followed by the alternative fact that the shape operator $\mathcal{S}_4 = \left(\begin{array}{ccc}
			0 & 0 & 0\\ 
			0 & 0 & 0 \\ 
			0 & 0 & 0\\ 
		\end{array}\right)$ attains the desired form.
\end{example}

Now, we give an example of a Riemannian map for which the inequality (\ref{eq-GSSF-(1)}) holds strictly.

\begin{example}
    Let $M=\left\{(x_1,x_2,\dots,x_6)\in\mathbb{R}^{6}\right\}$ and $N=\left\{(y_1,y_2,\dots, y_6)\in\mathbb{R}^{6}:y_3>-1,\;y_5>-1\right\}$ be two Riemannian manifolds with Riemannian metrics $g_1= \sum \limits_{i=1}^{6}dx_i^2$ and $g_2=dy_1^2+(1+y_3)^2dy_2^2+dy_3^2 +(1+y_5)^2dy_4^2+dy_5^2+dy_6^2$, respectively. Define a map $ \pi: \left( M, g_{1}\right) \to \left( N,g_{2}\right) $ by
    \[
    \pi \left( x_{1},x_{2},x_{3},x_{4},x_{5},x_{6}\right) =\left(x_1, x_2,0,x_{4},0,x_{6}\right).
    \]
    Then we have
    \begin{equation*}
        \mathcal{H} = \operatorname{span} \{h_1 = e_1, h_2 = e_2, h_3 = e_4, h_4 = e_6\}, \quad \mathcal{V} = \operatorname{span} \{e_3, e_5\},
    \end{equation*}
    \begin{equation*}
        \mathcal{R} =\operatorname{span} \left\{ \pi_{\ast } h_{1}=e_{1}^{\ast }, \pi_{\ast } h_{2} =e_{2}^{\ast }, \pi_{\ast}  h_{3} =e_{4}^{\ast }, \pi_{\ast} h_4 = e_6^{\ast} \right\}, \quad \mathcal{R}^{\perp} = \operatorname{span} \left\{V_5 =e_{3}^{\ast }, V_6 = e_5^{\ast} \right\},
    \end{equation*}
    where $\{e_i = \frac{\partial}{\partial x_i}\}_{i=1}^{6}$ and $\{e_j^{\ast} = \frac{\partial}{\partial y_j}\}_{j=1}^{6}$ are bases of $T_p M$ and $T_{\pi (p)} N$, respectively. We observe that $(\ref{riemannian_map})$ holds for all $h_{i}, h_{j} \in \Gamma (\mathcal{H})$ as $y_3 = y_5 = 0$ under the defined map. Thus, $\pi$ is a Riemannian map with $r=4$. Now, we compute the following non-zero Christoffel symbols for $g_2$ as
    \begin{eqnarray*}
    \Gamma_{y_2 y_3}^{y_2} =\Gamma_{y_3 y_2}^{y_2} = \frac{1}{1+y_3}, \quad  \Gamma_{y_2 y_2}^{y_3}= -(1 + y_3), \quad \Gamma_{y_4 y_5}^{y_4} = \Gamma_{y_5 y_4}^{y_4} =\frac{1}{1+y_5}, \quad \Gamma_{y_4 y_4}^{y_5}= -(1 + y_5).
    \end{eqnarray*}
Then, using $(\ref{eqn_sff})$ and $y_3=y_5=0$, we obtain
\[
\left( \nabla \pi_{\ast }\right) \left( h_{1}, h_{1}\right) =0, \left( \nabla \pi_{\ast }\right) \left( h_{2}, h_{2}\right) = - e_3^{\ast}, \left( \nabla \pi_{\ast }\right) \left( h_{3}, h_{3}\right) = - e_5^{\ast}, \left( \nabla \pi_{\ast }\right) \left( h_{4}, h_{4}\right) =0,
\]
and
\[
\left( \nabla \pi_{\ast }\right) \left( h_{i}, h_{j}\right) =0,\ {\rm for\ all\ } h_{i} \neq h_{j} \in \Gamma(\mathcal{H}).
\]
Consequently, only non-zero components are $\beta_{22}^5 = -1$ and $\beta_{33}^6  = -1$. Further, we compute
\begin{equation*}
    \|\tau^{\mathcal{H}} \|^2 = 2,
\end{equation*}
\begin{equation*}
    K^{\mathcal{R}} (\mathbb{P}) = g_2(R^N(\pi_\ast h_1, \pi_\ast h_2) \pi_\ast h_2, \pi_\ast h_1) = g_2(R^N(e_1^*, e_2^*)e_2^*,e_1^*) = 0,
\end{equation*}
\begin{equation*}
    2\rho^{\mathcal{R}} = \sum\limits_{i, j= 1}^{4} g_2(R^N(\pi_\ast h_i, \pi_\ast h_j) \pi_\ast h_j, \pi_\ast h_i) = 0,
\end{equation*}
\begin{equation*}
    K^{\mathcal{H}} (\mathbb{P}) = g_1(R^M(h_1, h_2) h_2, h_1) = g_1(R^M(e_1, e_2)e_2,e_1) = 0,
\end{equation*}
and
\begin{equation*}
    2\rho^{\mathcal{H}} = \sum\limits_{i, j= 1}^{4} g_1(R^M(h_i, h_j) h_j, h_i) =0.
\end{equation*}
Thus, the inequality in $(\ref{eq-GSSF-(1)})$ holds strictly as $0> -\frac{2}{3}$. This is also followed by the alternative fact that the shape operator \begin{equation*}
    \mathcal{S}_5 = \left(\begin{array}{cccc}
			0 & 0 & 0 &0\\ 
			0 & -1 & 0 &0\\ 
			0 & 0 & 0 &0 \\
            0 & 0 & 0 &0 \\
		\end{array}\right) \quad \text{and} \quad \mathcal{S}_6 = \left(\begin{array}{cccc}
			0 & 0 & 0 &0\\ 
			0 & 0 & 0 &0\\ 
			0 & 0 & -1 &0 \\
            0 & 0 & 0 &0 \\
		\end{array}\right)
\end{equation*} does not attain the desired equality form.
\end{example}
	
	\section{Utilization of general Chen's first inequality}\label{sec_utilization}
	The main use of general Chen's first inequality is to establish Chen's first inequality with other target space forms. In this section, we prove Theorem \ref{Theorem_GCSF_1} (generalized complex space form as the target) and Theorem \ref{theorem_GSSF} (generalized Sasakian space form as the target) using $(\ref{eq-GSSF-(1)})$.\\
	
	First, we prove Theorem \ref{Theorem_GCSF_1} as follows:
	
	\begin{proof}[Proof of Theorem $\ref{Theorem_GCSF_1}$]
		By (\ref{curvature_tensor_gcsf}) and (\ref{decomposition_with_J}), we have 
		\begin{eqnarray}\label{eq-GCSF-(4)}
			&&g_2(R^{N}(\pi_{\ast }X,\pi_{\ast }Y)\pi_{\ast }Z,\pi_{\ast }W) = f_{1}\{ g_{2}(\pi_{\ast}Y,\pi_{\ast }Z)g_{2}(\pi_{\ast }X,\pi_{\ast }W) \nonumber \\&& -g_{2}(\pi_{\ast }X,\pi_{\ast }Z)g_{2}(\pi_{\ast}Y,\pi_{\ast }W)\}  \nonumber \\&&+f_{2}\left\{ g_{2}\left( \pi_{\ast }X, \mathcal{P} \pi_\ast Z\right) g_{2}\left( \mathcal{P} \pi_{\ast }Y,\pi_{\ast}W\right) -g_{2}\left( \pi_{\ast }Y, \mathcal{P} \pi_{\ast}Z\right) g_{2}\left( \mathcal{P} \pi_{\ast }X,\pi_{\ast}W\right) \right.  \nonumber \\&&\left. +2g_{2}\left( \pi_{\ast }X, \mathcal{P} \pi_{\ast }Y\right) g_{2}\left( \mathcal{P} \pi_{\ast }Z,\pi_{\ast}W\right) \right\}.
		\end{eqnarray}
		Since $\pi$ is a Riemannian map, (\ref{riemannian_map}) and (\ref{eq-GCSF-(4)}) yield,
		\begin{eqnarray*}
			&&g_2(R^{N}(\pi_{\ast }X,\pi_{\ast }Y)\pi_{\ast }Z,\pi_{\ast }W) =f_{1}\{g_{1}(Y,Z)g_{1}(X,W)-g_{1}(X,Z)g_{1}(Y,W)\} \\&&+f_{2}\left\{ g_{2}\left( \pi_{\ast }X, \mathcal{P} \pi_{\ast }Z\right) g_{2}\left( \mathcal{P} \pi_{\ast}Y,\pi_{\ast }W\right) -g_{2}\left( \pi_{\ast }Y, \mathcal{P} \pi_{\ast }Z\right) g_{2}\left( \mathcal{P} \pi_{\ast}X,\pi_{\ast }W\right) \right. \\&&\left. +2g_{2}\left( \pi_{\ast }X, \mathcal{P} \pi_{\ast }Y\right) g_{2}\left( \mathcal{P} \pi_{\ast }Z,\pi_{\ast}W\right) \right\}.
		\end{eqnarray*}
		Contracting the above equation, we obtain
		\begin{eqnarray*}
			\sum \limits_{i,j=1}^{r} g_2(R^{N}(\pi_{\ast }h_{i},\pi_{\ast }h_{j})\pi_{\ast }h_{j},\pi_{\ast }h_{i}) =&\sum \limits_{i,j=1}^{r}f_{1} \{ g_{1}(h_{j},h_{j})g_{1}(h_{i},h_{i})-g_{1}(h_{i},h_{j})g_{1}(h_{j},h_{i})\} \\&+\sum \limits_{i,j=1}^{r}f_{2}\left\{ g_{2}\left( \pi_{\ast }h_{i}, \mathcal{P} \pi_{\ast }h_{j}\right) g_{2}\left( \mathcal{P} \pi_{\ast }h_{j},\pi_{\ast }h_{i}\right) \right. \\&\left.-g_{2}\left( \pi_{\ast }h_{j}, \mathcal{P} \pi_{\ast}h_{j}\right) g_{2}\left( \mathcal{P} \pi_{\ast }h_{i},\pi_{\ast }h_{i}\right) \right. \\&\left. +2g_{2}\left( \pi_{\ast }h_{i}, \mathcal{P} \pi_{\ast }h_{j}\right) g_{2}\left(\mathcal{P} \pi_{\ast }h_{j},\pi_{\ast }h_{i}\right) \right\},
		\end{eqnarray*}
		where $\left\{ h_{i}\right\}_{i=1}^{r}$ (resp. $\left\{ \pi_{\ast}h_{i}\right\}_{i=1}^{r}$) is an orthonormal basis of $\mathcal{H}$ (resp. $\mathcal{R}$). This gives
		\begin{equation}\label{eq-GCSF-(5)}
			2 \rho^{\mathcal{R}} = r \left(r-1\right) f_{1} + 3 f_{2} \|\mathcal{P}\|^{2}.
		\end{equation}
		In particular, if $\pi_{\ast }X=\pi_{\ast }W=\pi_{\ast }h_{1}$ and $\pi_{\ast}Y = \pi_{\ast }Z = \pi_{\ast }h_{2}$ in (\ref{eq-GCSF-(4)}), then we obtain
		\begin{equation}\label{eq-GCSF-(6)}
			K^{\mathcal{R}}\left( \mathbb{P} \right) =f_{1}+3f_{2} (g_{2}\left( \pi_{\ast }h_{1}, \mathcal{P} \pi_{\ast }h_{2}\right))^{2}.
		\end{equation}
        Now, we discuss the following cases:
        \begin{enumerate}[$(1)$]
            \item If $f_1$ and $f_2$ are constant, then we have $3 \leq r < n$. Therefore, using (\ref{eq-GCSF-(5)}) and (\ref{eq-GCSF-(6)}) in (\ref{eq-GSSF-(1)}), we get (\ref{eq-GCSF-(1)}). Hence, the proof of $(1)$ follows.
            \item If $f_1$ and $f_2$ are non-constant, then due to $\cite{Olszak}$ and Remark $\ref{rem_on_rank}$, we must have $n = 4$ and $r=3$. In this case, the range space will be a hypersurface of the target space \cite{Roth_Upadhyay}. Therefore, using (\ref{eq-GCSF-(5)}) and (\ref{eq-GCSF-(6)}) in (\ref{eq-GSSF-(1)}) with substituting $r=3$, we get (\ref{const_case_gcsf}). Hence, the proof of $(2)$ follows.
        \end{enumerate}
        Further, the equality cases are straightforward due to the equality case of (\ref{eq-GSSF-(1)}).
	\end{proof}
	
Along similar lines, we now prove Theorem \ref{theorem_GSSF} as follows:
	
	\begin{proof}[Proof of Theorem $\ref{theorem_GSSF}$]
		If $\xi \in \Gamma(\mathcal{R})$, then by (\ref{curvature_tensor_gssf}) and (\ref{decomposition_with_J}), we have
		\begin{eqnarray}\label{eq-GNSSF-(4)}
			&&g_2(R^N \left( \pi_{\ast }X,\pi_{\ast }Y)\pi_{\ast }Z,\pi_{\ast }W\right)\nonumber\\&&=f_{1}\left\{ g_{2}\left( \pi_{\ast }Y,\pi_{\ast }Z\right) g_{2}\left( \pi_{\ast}X,\pi_{\ast }W\right) -g_{2}\left( \pi_{\ast }X,\pi_{\ast }Z\right) g_{2}\left(\pi_{\ast }Y,\pi_{\ast }W\right) \right\} \nonumber\\&&+f_{2}\left\{ g_{2}\left( \pi_{\ast }X, \mathcal{P} \pi_{\ast }Z\right) g_{2}\left(\mathcal{P} \pi_{\ast }Y, \pi_{\ast }W\right) -g_{2}\left( \mathcal{P} \pi_{\ast }X, \pi_{\ast}W\right) g_{2}\left( \pi_{\ast }Y,\mathcal{P} \pi_{\ast }Z\right) \right.   \nonumber\\&&+\left. 2g_{2}\left( \pi_{\ast }X,\mathcal{P} \pi_{\ast }Y\right) g_{2}\left( \mathcal{P} \pi_{\ast}Z, \pi_{\ast }W\right) \right\}   \nonumber \\&& +f_{3}\left\{ g_{2}\left( \pi_{\ast }Y,\pi_{\ast }W\right) \eta \left( \pi_{\ast}X\right) \eta \left( \pi_{\ast }Z\right) -g_{2}\left( \pi_{\ast }X,\pi_{\ast}W\right) \eta \left( \pi_{\ast }Y\right) \eta \left( \pi_{\ast }Z\right)\right.   \nonumber \\&&\left. +g_{2}\left( \pi_{\ast }X,\pi_{\ast }Z\right) \eta \left( \pi_{\ast}Y\right) \eta \left( \pi_{\ast }W\right) -g_{2}\left( \pi_{\ast }Y,\pi_{\ast}Z\right) \eta \left( \pi_{\ast }X\right) \eta \left( \pi_{\ast }W\right)\right\}.
		\end{eqnarray}
		Since $\pi$ is a Riemannian map, (\ref{riemannian_map}) and (\ref{eq-GNSSF-(4)}) yield,
		\begin{eqnarray*}
			&&g_2(R^N\left( \pi_{\ast }X,\pi_{\ast }Y)\pi_{\ast }Z,\pi_{\ast }W\right)\nonumber\\&&=f_{1}\left\{ g_{1}\left( Y,Z\right) g_{1}\left( X,W\right) -g_{1}\left(X,Z\right) g_{1}\left( Y,W\right) \right\}  \\&&+f_{2}\left\{ g_{2}\left( \pi_{\ast }X,\mathcal{P} \pi_{\ast }Z\right) g_{2}\left(\mathcal{P} \pi_{\ast }Y,\pi_{\ast }W\right) -g_{2}\left( \pi_{\ast }Y,\mathcal{P} \pi_{\ast}Z\right) g_{2}\left( \mathcal{P} \pi_{\ast }X,\pi_{\ast }W\right) \right.  \\&&+\left. 2g_{2}\left( \pi_{\ast }X, \mathcal{P} \pi_{\ast }Y\right) g_{2}\left( \mathcal{P} \pi_{\ast} Z, \pi_{\ast }W\right) \right\}  \\&&+f_{3}\left\{ \eta \left( \pi_{\ast }X\right) \eta \left( \pi_{\ast }Z\right)g_{2}\left( \pi_{\ast }Y,\pi_{\ast }W\right) -\eta \left( \pi_{\ast }Y\right) \eta\left( \pi_{\ast }Z\right) g_{2}\left( \pi_{\ast }X,\pi_{\ast }W\right) \right.  \\&&\left. +g_{2}\left( \pi_{\ast }X,\pi_{\ast }Z\right) \eta \left( \pi_{\ast}Y\right) \eta \left( \pi_{\ast }W\right) -g_{2}\left( \pi_{\ast }Y,\pi_{\ast}Z\right) \eta \left( \pi_{\ast }X\right) \eta \left( \pi_{\ast }W\right)\right\}. 
		\end{eqnarray*}
		Contracting the aforementioned equation, we get
		\begin{eqnarray*}
			&&\sum \limits_{i,j=1}^{r} g_2(R^N\left( \pi_{\ast }h_{i},\pi_{\ast }h_{j})\pi_{\ast}h_{j},\pi_{\ast }h_{i}\right)\\&&=\sum \limits_{i,j=1}^{r}f_{1}\left\{ g_{1}\left( h_{j},h_{j}\right)g_{1}\left( h_{i},h_{i}\right) -g_{1}\left( h_{i},h_{j}\right) g_{1}\left(h_{j},h_{i}\right) \right\}  \\&&+\sum \limits_{i,j=1}^{r}f_{2}\left\{ g_{2}\left( \pi_{\ast }h_{i},\mathcal{P} \pi_{\ast }h_{j}\right)g_{2}\left( \mathcal{P} \pi_{\ast }h_{j},\pi_{\ast }h_{i}\right) \right.  \\&& \left. -g_{2}\left( \pi_{\ast}h_{j},\mathcal{P} \pi_{\ast }h_{j}\right) g_{2}\left( \mathcal{P} \pi_{\ast }h_{i},\pi_{\ast}h_{i}\right) + 2g_{2}\left( \pi_{\ast }h_{i},\mathcal{P} \pi_{\ast }h_{j}\right) g_{2}\left(\mathcal{P} \pi_{\ast }h_{j},\pi_{\ast }h_{i}\right) \right\}  \\&&+\sum \limits_{i,j=1}^{r}f_{3}\left\{ \eta \left( \pi_{\ast }h_{i}\right) \eta \left( \pi_{\ast}h_{j}\right) g_{2}\left( \pi_{\ast }h_{j},\pi_{\ast }h_{i}\right) -\eta \left(\pi_{\ast }h_{j}\right) \eta \left( \pi_{\ast }h_{j}\right) g_{2}\left( \pi_{\ast}h_{i},\pi_{\ast }h_{i}\right) \right.  \\&&\left. +g_{2}\left( \pi_{\ast }h_{i},\pi_{\ast }h_{j}\right) \eta \left(\pi_{\ast }h_{j}\right) \eta \left( \pi_{\ast }h_{i}\right) -g_{2}\left( \pi_{\ast}h_{j},\pi_{\ast }h_{j}\right) \eta \left( \pi_{\ast }h_{i}\right) \eta \left(\pi_{\ast }h_{i}\right) \right\},
		\end{eqnarray*}
		where $\left\{ h_{i}\right\}_{i=1}^{r}$ (resp. $\left\{ \pi_{\ast}h_{i}\right\}_{i=1}^{r}$) is an orthonormal basis of $\mathcal{H}$ (resp. $\mathcal{R}$). This gives
		\begin{equation}\label{eq-GNSSF-(5)}
			2\rho^{\mathcal{R}} = r\left( r-1\right) f_1 + 3 f_{2} \|\mathcal{P}\|^{2} -2 (r-1) f_3.
		\end{equation}
		In particular, if $\pi_{\ast }X=\pi_{\ast }W=\pi_{\ast}h_{1}$ and $\pi_{\ast}Y=\pi_{\ast }Z=\pi_{\ast }h_{2}$ in (\ref{eq-GNSSF-(4)}), then we obtain
		\begin{equation}\label{eq-GNSSF-(6)}
			K^{\mathcal{R}}\left( \mathbb{P} \right) =f_{1}+3f_{2} \left(g_{2}\left(\pi_{\ast }h_{1},\mathcal{P} \pi_{\ast }h_{2}\right)\right)^{2} - f_3 \left\{\left(\eta(\pi_\ast h_1)\right)^2 + \left(\eta(\pi_\ast h_2)\right)^2\right\}.
		\end{equation}
		Using (\ref{eq-GNSSF-(5)}) and (\ref{eq-GNSSF-(6)}) in (\ref{eq-GSSF-(1)}), we get (\ref{eq-GNSSF-(1)}).\\
		
		Further, we observe that if $\xi \in \Gamma(\mathcal{R}^\perp)$, then the coefficients of $f_3$ in $(\ref{eq-GNSSF-(5)})$ and $(\ref{eq-GNSSF-(6)})$ vanish. Hence (\ref{CFI_when_xi_in_R_perp}) follows.\\
        
       Furthermore, the equality cases are straightforward due to the equality case of (\ref{eq-GSSF-(1)}).
	\end{proof}
	
	\section{Some important consequences}\label{sec_corollary}
	In the previous section, we established Chen's first inequalities when target spaces admit generalized complex and generalized Sasakian structures. The umbrellas of generalized complex and generalized Sasakian space forms include the following specific space forms, depending on the following values of $f_1, f_2$ and $f_3$ \cite{ABC, Vanhecke, Tricerri}:
	\begin{center}
		\begin{tabular}{|c|c|c|c|}
			\hline
			\textbf{space form}& \textbf{$f_1$} & \textbf{$f_2$}\\
			\hline real  & $c$ & $0$\\
			complex & $c$ & $c$\\
			real K\"ahler &$c+3\alpha$& $c-\alpha$\\
			\hline
		\end{tabular}
		\begin{tabular}{|c|c|c|}
			\hline
			\textbf{space form}& \textbf{$f_1$} & \textbf{$f_2=f_3$} \\
			\hline Sasakian  & $c+3$ & $c-1$\\
			Kenmotsu  & $c-3$ & $c+1$\\
			cosymplectic  & $c$ & $c$\\
			almost $C(\alpha)$  &$c+3\alpha^2$& $c-\alpha^2$\\
			\hline
		\end{tabular}
	\end{center}
	Thus, these space forms unify various structures and facilitate studying Chen's first inequalities across them. In this section, we apply Chen's first inequalities (Theorems \ref{Theorem_GCSF_1} and \ref{theorem_GSSF}) with direct substitution to explore a framework that adapts several existing and dominant cases.\\ 
    
    Accordingly, we have the following corollary of Theorem \ref{Theorem_GCSF_1}:
	\begin{corollary}\label{Corollary_GCSF}
		Let $\pi:\left( M^{m},g_{1}\right) \to \left(N^n, g_{2}\right)$ be a Riemannian map. Then the following statements are valid:
		\begin{enumerate}[$(1)$]
			\item If $N\left(c\right)$ is a real space form, then
			\begin{equation*}
				K^{\mathcal{H}}\left( \mathbb{P} \right) \geq \rho^{\mathcal{H}}-\frac{(r-2)}{2}\left( \left( r+1\right)c+\frac{1}{(r-1)}\|\tau^{\mathcal{H}}\|^{2}\right).
			\end{equation*}
			
			\item If $N(4c)$ is a complex space form, then
			\begin{equation*}
				K^{\mathcal{H}}\left( \mathbb{P} \right) \geq \frac{1}{2}\left\{ 2\rho^{\mathcal{H}}-\frac{\left(r-2\right) }{\left( r-1\right) }\|\tau^{\mathcal{H}}\|^{2}-c\left( r^{2}-r-2\right) -3c\left(\| \mathcal{P} \|^{2}-2 \Theta(\mathbb{P})\right) \right\}.
			\end{equation*}
			
			\item If $N^{n=4}(4c)$ is a real K\"ahler space form, then
			\begin{align*}
				&K^{\mathcal{H}}\left( \mathbb{P} \right) \geq  \frac{1}{2}\left\{ 2\rho^{\mathcal{H}}-\frac{1}{2} \|\tau^{\mathcal{H}}\|^{2}-4(c+3\alpha) -3(c-\alpha)\left( \| \mathcal{P} \|^{2}-2 \Theta(\mathbb{P}) \right) \right\}.
			\end{align*}
		\end{enumerate}
		In addition, the equality conditions are identical to those of $(\ref{eq-GSSF-(1)})$.
	\end{corollary}
    
	\begin{remark}
		Observe that since we denote scalar curvature by $2 \rho^{\mathcal{H}}$ instead of $\rho^{\mathcal{H}}$, the inequalities in items $(1)$ and $(2)$ of Corollary $\ref{Corollary_GCSF}$ are the same as the inequalities given in $\cite[~{\rm Theorem} ~3.1]{Sahin_2016}$ and $\cite[~{\rm Theorem} ~3.1]{Meena}$, respectively.
	\end{remark}

    In addition, we have the following corollary of Theorem \ref{theorem_GSSF}:
    \begin{corollary}\label{Corollary_1_GSSF}
		Let $\pi:\left( M^{m},g_{1}\right) \to \left( N^{n}, g_{2}, \phi, \xi, \eta\right) $ be a Riemannian map. Then the following statements are valid:
		\begin{enumerate}[$(1)$]
			\item Suppose $N(4c)$ is a Sasakian space form, then
			
			\begin{equation*}
				K^{\mathcal{H}}\left( \mathbb{P} \right) - \rho^{\mathcal{H}} +\frac{\left(r-2\right) }{2\left( r-1\right) }\|\tau^{\mathcal{H}}\|^{2}\geq 
				\left\{\begin{array}{ll}
					-\frac{r((r-3)c + 3r -1)}{2} +4- \frac{(c-1)}{2} \left(3\|\mathcal{P}\|^{2}-6~ \Theta(\mathbb{P}) + 2 \Phi(\mathbb{P})\right), ~{\rm if}~ \xi \in \Gamma(\mathcal{R});\\- (r^2-r-2)(c+3) -3 (c-1)\left( \|\mathcal{P}\|^{2}-2\Theta(\mathbb{P}) \right), ~{\rm if}~ \xi \in \Gamma(\mathcal{R}^\perp).	
				\end{array} \right.
			\end{equation*}
			
			\item Suppose $N(4c)$ is a Kenmotsu space form, then
			
			\begin{equation*}
				K^{\mathcal{H}}\left( \mathbb{P} \right) - \rho^{\mathcal{H}} +\frac{\left(r-2\right) }{2\left( r-1\right) }\|\tau^{\mathcal{H}}\|^{2}\geq 
				\left\{\begin{array}{ll}
					-\frac{r((r-3)c-3r+1)}{2} -4-\frac{(c+1)}{2} \left(3\|\mathcal{P}\|^{2}-6~ \Theta(\mathbb{P}) + 2 \Phi(\mathbb{P})\right), ~{\rm if}~ \xi \in \Gamma(\mathcal{R});\\- (r^2-r-2)(c-3) -3 (c+1)\left( \|\mathcal{P}\|^{2}-2\Theta(\mathbb{P}) \right), ~{\rm if}~ \xi \in \Gamma(\mathcal{R}^\perp).	
				\end{array} \right.
			\end{equation*}
			
			\item Suppose $N(4c)$ is a cosymplectic space form, then
			\begin{equation*}
				K^{\mathcal{H}}\left( \mathbb{P} \right) - \rho^{\mathcal{H}} +\frac{\left(r-2\right) }{2\left( r-1\right) }\|\tau^{\mathcal{H}}\|^{2}\geq 
				\left\{\begin{array}{ll}
					-\frac{c}{2}r\left( r-3\right) -\frac{c}{2} \left(3\|\mathcal{P}\|^{2}-6~ \Theta(\mathbb{P}) + 2 \Phi(\mathbb{P})\right), ~{\rm if}~ \xi \in \Gamma(\mathcal{R});\\- (r^2-r-2)c -3 c\left( \|\mathcal{P}\|^{2}-2\Theta(\mathbb{P}) \right), ~{\rm if}~ \xi \in \Gamma(\mathcal{R}^\perp).	
				\end{array} \right.
			\end{equation*}
			
			\item Suppose $N(4c)$ is an almost $C(\alpha)$ space form, then
			\begin{align*}
				K^{\mathcal{H}}\left( \mathbb{P} \right) &- \rho^{\mathcal{H}} +\frac{\left(r-2\right) }{2\left( r-1\right) }\|\tau^{\mathcal{H}}\|^{2} \\& \geq 
				\left\{\begin{array}{ll}
					- \frac{r((r-3)c +\alpha^2 (3r-1))}{2} + 4 \alpha^2 -\frac{(c-\alpha^2)}{2} \left(3\|\mathcal{P}\|^{2}-6~ \Theta(\mathbb{P}) + 2 \Phi(\mathbb{P})\right), ~{\rm if}~ \xi \in \Gamma(\mathcal{R});\\- (r^2-r-2)(c+3\alpha^2) -3 (c-\alpha^2)\left( \|\mathcal{P}\|^{2}-2\Theta(\mathbb{P}) \right), ~{\rm if}~ \xi \in \Gamma(\mathcal{R}^\perp).	
				\end{array} \right.
			\end{align*}
		\end{enumerate}
		In addition, the equality conditions are identical to those of $(\ref{eq-GSSF-(1)})$.
	\end{corollary}
	
	\section{Estimations of $\delta$-invariants}\label{sec_deltainv}
	In this section, we present geometric applications by estimating $\delta$-invariants. By \cite{ACM, Chen_RIMA, CDVV_1996}, we know that \textit{$\delta$-invariants} are defined as
	\begin{equation*}
		\delta^{\mathcal{H}}= \rho^{\mathcal{H}} - \inf\{K^{\mathcal{H}}(\mathbb{P}) : \text{plane sections}~ \mathbb{P} \subset \mathcal{H}_p\}.
	\end{equation*}
	Before estimating them, we list the following statements for ease of presentation.\\
	
	\noindent There exist orthonormal bases $\{\pi_{\ast}h_{i}\}_{i=1}^{r}$ and $\{V_{k}\}_{k=r+1}^{n}$ of $\mathcal{R}$ and $\mathcal{R}^\perp$, respectively, such that
	\begin{enumerate}[$A:$]
		\item leaves of $\mathcal{R}$ are even-dimensional invariant almost Hermitian submanifolds of the target space,\label{S1}
		
        \item $\pi_\ast h_1$ and $\mathcal{P} \pi_\ast h_2$ are orthogonal,\label{S2}
		
        \item the shape operator satisfies $(\ref{eq-GSSF-(2)})$,\label{S3}
		
        \item $K^{{\mathcal{H}}}(h_{1} \wedge h_{2})=\inf \{K^{{\mathcal{H}}}(\mathbb{P} ):{\rm \ }${\rm plane\ sections}$\ \mathbb{P} \subset \mathcal{H}_{p}\}$,\label{S4}
		
        \item $\pi_\ast h_1, \pi_\ast h_2, \mathcal{P} \pi_\ast h_3, \mathcal{P} \pi_\ast h_4, \dots, \mathcal{P} \pi_\ast h_r$ are orthogonal and the subspace spanned by $\{\pi_{\ast }h_{i}\}_{i=3}^{r}$ is totally real, equivalently, $\mathcal{P} \pi_\ast h_i = 0$ for all $3 \leq i \leq r$. In particular, in the sequel when $r=3$, it can be interpret as either $\pi_\ast h_1, \pi_\ast h_2$, and $\mathcal{P} \pi_\ast h_3$ are orthogonal, or $\mathcal{P}\pi_\ast h_3 = 0$,\label{S5}
		
        \item $\xi, \pi_\ast h_1$, and $\pi_\ast h_2$ are orthonormal,\label{S6}
		
		\item leaves of $\mathcal{R}$ are odd-dimensional invariant almost contact submanifolds of the target space,\label{S7}
		
        \item $ \xi, \pi_\ast h_3, \pi_\ast h_4, \dots$, $\pi_\ast h_r $ are orthonormal,\label{S8}
	\end{enumerate}
	
	First, we provide $\delta$-invariant estimates for Riemannian maps to generalized complex space forms. By \cite{Olszak, Tricerri}, we know that every generalized complex space form of dimension $n>4$ with $f_2 \neq 0$ reduces to a complex space form. Since the $\delta$-invariants already exist for complex space forms \cite{Meena}, we restrict ourselves to $n=4$ in the following theorem. In this case, the $\delta$-invariant estimates do not depend on the constancy of $f_1$ and $f_2$.
	\begin{theorem}\label{theorem_f2geq0}
		Let $\pi:(M^{m},g_{1})\to (N^{n=4}(f_{1},f_{2}),g_{2}, J)$ be a Riemannian map from a Riemannian manifold to a generalized complex space form.
		\begin{enumerate}[$(1)$]
			\item Suppose $f_{2}> 0$, then
			\begin{equation}\label{eq-CGCSF-(1)}
				\delta^{{\mathcal{H}}} < \frac{1}{2} \left\{ \frac{1}{2}\| \tau^{{\mathcal{H}}}\|^2 +4 f_{1} + 9 f_{2}\right\}.
			\end{equation}
			
			\item Suppose $f_{2}\leq 0$, then
			\begin{equation}\label{eq-CGCSF-(2)}
				\delta^{{\mathcal{H}}}\leq \frac{\| \tau^{{\mathcal{H}}}\|^2}{4} +2f_{1}.
			\end{equation}
			The equality in $(\ref{eq-CGCSF-(2)})$ holds if and only if the statements
            $$\left\{\begin{array}{ll}
			\text{\ref{S3}, \ref{S4}, and \ref{S5} follow}, & \text{when }~ f_2 < 0; \\
			\text{\ref{S3} and \ref{S4} follow}, & \text{when }~ f_2 = 0.
		\end{array} \right.$$
		\end{enumerate}
	\end{theorem}
	
	\begin{proof}
    As $n=4$, due to Remark \ref{rem_on_rank}, we have $r=3$. Consequently, in this case, the inequalities mentioned in (\ref{eq-GCSF-(1)}) and (\ref{const_case_gcsf}) coincide with the following inequality:
    \begin{equation}\label{new_eq_delta_H_gcsf}
				\delta^{{\mathcal{H}}}\leq \frac{1}{2}\left\{ \frac{1}{2}\| \tau^{{\mathcal{H}}}\|^2 +4 f_{1}+3f_{2}\left( \| \mathcal{P} \|^{2}-2\Theta (\mathbb{P})\right) \right\}.
			\end{equation}
            In addition, the equality in the above inequality holds if and only if the statements $\ref{S3}$ and $\ref{S4}$ hold. We now prove the required estimates one by one.
		\begin{enumerate}[$(1)$]
			\item To maximize the right-hand side of (\ref{new_eq_delta_H_gcsf}), we have to maximize the term $f_2 \left(\|\mathcal{P}\|^{2}-2\Theta (\mathbb{P} )\right)$ for some fixed $f_1$. Since $f_{2}> 0$, we must maximize $\left(\|\mathcal{P}\|^{2}-2\Theta (\mathbb{P} )\right)$. Obviously, it attains its maximum when $\|\mathcal{P}\|^2$ is maximized and $\Theta(\mathbb{P})$ is minimized. For their compatibility, we have the following key facts:
            \begin{itemize}
                \item We can write $\mathcal{P} \pi_\ast h_j = \sum\limits_{i=1}^3 g_2(\mathcal{P} \pi_\ast h_j, \pi_\ast h_i)\pi_\ast h_i$; therefore $\|\mathcal{P} \pi_\ast h_j \|^2 = \sum \limits_{i=1}^3 (g_2(\pi_\ast h_i, \mathcal{P} \pi_\ast h_j))^2$. Thus,
                \begin{equation}\label{norm_P_square}
                    \sum \limits_{j=1}^{3} \|\mathcal{P} \pi_\ast h_j \|^2 = \| \mathcal{P} \|^2.
                \end{equation}
                We also have
                \begin{equation*}
                    \| J \pi_\ast h_j \|^2 = g_2(J \pi_\ast h_j, J \pi_\ast h_j) = -g_2(J^2 \pi_\ast h_j, \pi_\ast h_j) = 1 \quad [\text{since $J^2 = -I$}]
                \end{equation*}
                and 
                \begin{align*}
                    \| J \pi_\ast h_j \|^2 &= g_2(J \pi_\ast h_j, J \pi_\ast h_j) \\&= g_2(\mathcal{P} \pi_\ast h_j + \mathcal{Q} \pi_\ast h_j, \mathcal{P} \pi_\ast h_j + \mathcal{Q} \pi_\ast h_j) \quad [\text{using (\ref{decomposition_with_J})}] \\& = \| \mathcal{P} \pi_\ast h_j \|^2 + \| \mathcal{Q} \pi_\ast h_j \|^2.
                \end{align*}
                The above expressions mean $\| \mathcal{P} \pi_\ast h_j \|^2 + \| \mathcal{Q} \pi_\ast h_j \|^2 = 1$; hence $\| \mathcal{P} \pi_\ast h_j \|^2 \leq 1$. Consequently, (\ref{norm_P_square}) gives $\| \mathcal{P}\|^2 \leq 3$. Clearly, the maximum value of $\| \mathcal{P}\|^2$ may reaches $3$. 

                \item By definition, we have $\Theta(\mathbb{P}) =\left(g_{2}( \pi_{\ast }h_{1}, \mathcal{P} \pi_{\ast }h_{2})\right)^{2} \geq 0$. Clearly, the minimum value of $\Theta(\mathbb{P})$ reaches $0$.
            \end{itemize}
            Thus, using the above key facts, we obtain
            \begin{equation*}
				\delta^{{\mathcal{H}}} \leq \frac{1}{2} \left\{ \frac{1}{2}\| \tau^{{\mathcal{H}}}\|^2 +4 f_{1} + 9 f_{2}\right\}.
			\end{equation*}
            In addition, equality in the above inequality requires $\| \mathcal{P}\|^2 =3$ and $\Theta(\mathbb{P}) = 0$. For such a case, one can observe that the condition $\| \mathcal{P} \|^{2}=3$ together with $\| \mathcal{P} \pi_\ast h_j \|^2 + \| \mathcal{Q} \pi_\ast h_j \|^2 = 1$ implies that $\mathcal{Q} \pi_\ast h_j = 0$ for each $1 \leq j \leq 3$. Then (\ref{decomposition_with_J}) yields
                \begin{equation*}
                    J \pi_\ast h_j = \mathcal{P} \pi_\ast h_j \quad \text{for all $1 \leq j \leq 3$}.
                \end{equation*}
                Equivalently, $J (\mathcal{R}) \subseteq (\mathcal{R})$; hence the leaves of $\mathcal{R}$ are $J$-invariant submanifolds of the target space \cite[Section 5]{MMT_conf_paper}. Further, by \cite[Theorem 5.4]{MMT_conf_paper}, they are almost Hermitian submanifolds, so $r$ must be even-dimensional. This is a contradiction as we have $r=3$. Therefore, the equality cannot occur. Thus, the strict inequality (\ref{eq-CGCSF-(1)}) follows.
			
			\item The Equation (\ref{new_eq_delta_H_gcsf}) can be simplified as
			\begin{align*}
				\delta^{{\mathcal{H}}} \leq&  \frac{\| \tau^{{\mathcal{H}}}\|^2}{4} + 2 f_{1} + \frac{3}{2} f_{2} \left\{(g_2(\pi_\ast h_1, \mathcal{P} \pi_\ast h_3))^2  + (g_2(\pi_\ast h_2, \mathcal{P} \pi_\ast h_3))^2\right\}.
			\end{align*}
			Since $f_2 \leq 0$, we must minimize the coefficient of $f_2$ in the aforementioned equation. Therefore, we can write (\ref{eq-CGCSF-(2)}) because the coefficient of $f_2$ reaches its minimum value $0$. When $f_2 <0$, the equality in (\ref{eq-CGCSF-(2)}) holds if and only if $(g_2(\pi_\ast h_1, \mathcal{P} \pi_\ast h_3))^2 + (g_2(\pi_\ast h_2, \mathcal{P} \pi_\ast h_3))^2  = 0$ and the statements $\ref{S3}$ and $\ref{S4}$ satisfy (due to equality in (\ref{new_eq_delta_H_gcsf})). It is straightforward that the condition $(g_2(\pi_\ast h_1, \mathcal{P} \pi_\ast h_3))^2 + (g_2(\pi_\ast h_2, \mathcal{P} \pi_\ast h_3))^2 = 0$ is equivalent to the statement $\ref{S5}$. It is also worth noting that if $f_2=0$, then the statement $\ref{S5}$ is not necessary for equality because in this case the term $$\frac{3}{2} f_{2} \left\{(g_2(\pi_\ast h_1, \mathcal{P} \pi_\ast h_3))^2  + (g_2(\pi_\ast h_2, \mathcal{P} \pi_\ast h_3))^2\right\}$$ vanishes identically.
		\end{enumerate}
	\end{proof}
	
	Using the fact that $\tau^\mathcal{H}=0$ for harmonic Riemannian maps, we have the following corollary:
	\begin{corollary}\label{cor_harmonic_gcsf}
		In particular, suppose $\pi:(M^{m},g_{1})\to (N^{n=4}(f_{1},f_{2}),g_{2},J)$ is a harmonic Riemannian map; then
		$$\delta^{{\mathcal{H}}} \leq 
		\left\{\begin{array}{ll}
			2 {f_1}, & \text{if }~ f_2 \leq 0;\\
			\frac{4 f_1 + 9 f_2}{2}, & \text{if }~ f_2 > 0.
		\end{array} \right.$$
	\end{corollary}

    \begin{remark}
    One can observe the following: 
    \begin{enumerate}[$(i)$]
        \item If $f_1 = f_2 = c > 0$, then $(\ref{eq-CGCSF-(1)})$ gives $\delta^{{\mathcal{H}}} < \frac{1}{2} \left\{ \frac{1}{2}\| \tau^{{\mathcal{H}}}\|^2 +13 c\right\}$. 
        
        \item If $f_1 = f_2 = c = 0$, then $(\ref{eq-CGCSF-(2)})$ gives $\delta^{{\mathcal{H}}} \leq \frac{1}{4}\| \tau^{{\mathcal{H}}}\|^2$.

        \item If $f_1 = f_2 = c < 0$, then $(\ref{eq-CGCSF-(2)})$ gives $\delta^{{\mathcal{H}}} \leq \frac{1}{4}\| \tau^{{\mathcal{H}}}\|^2 + 2c$.
    \end{enumerate}
    The same inequalities can also be obtained by particularly imposing $n=4$ and $r=3$ in $\cite[~{\rm Theorem}~ 4.2]{Meena}$, $\cite[~{\rm Theorem}~ 4.1]{Meena}$, and $\cite[~{\rm Theorem}~ 4.3]{Meena}$, respectively.
    \end{remark}
	
	Now, we provide $\delta$-invariant estimates for Riemannian maps to generalized Sasakian space forms. 
	
	\begin{theorem}\label{thm_gssf_1}
		Let $\pi:(M^{m},g_{1})\to (N^{n}(f_{1},f_{2},f_{3}),g_{2}, \phi, \xi, \eta)$ be a Riemannian map from a Riemannian manifold to a generalized Sasakian space form with Reeb vector field $\xi \in \Gamma(\mathcal{R})$. 		
		\begin{enumerate}[$(1)$]
			\item Suppose $f_{2}\leq 0$ and $f_{3} \leq 0$, then
			\begin{equation}\label{eq-CGSSF-(1)}
				\delta^{{\mathcal{H}}}\leq \frac{\left( r-2\right)}{2}\left\{ \frac{1 }{\left(r-1\right) }\| \tau^{{\mathcal{H}}}\|^2 +\left( r+1\right) f_{1} \right\} - f_{3}\left( r-1\right).
			\end{equation}
			 The equality in $(\ref{eq-CGSSF-(1)})$ holds if and only if the statements
            $$\left\{\begin{array}{ll}
			\text{\ref{S3}, \ref{S4}, \ref{S5}, and \ref{S6} follow}, & \text{when }~ f_2 < 0 ~\text{and}~ f_3 <0; \\
			\text{\ref{S3}, \ref{S4}, and \ref{S6} follow}, & \text{when }~ f_2 = 0 ~\text{and}~ f_3 <0; \\ \text{\ref{S3}, \ref{S4}, and \ref{S5} follow}, & \text{when }~ f_2 < 0 ~\text{and}~ f_3 =0; \\
			\text{\ref{S3} and \ref{S4} follow}, & \text{when }~ f_2 = 0 ~\text{and}~ f_3 =0.
		\end{array} \right.$$
			
			\item Suppose $f_{2}> 0$ and $f_{3} \leq 0$, then
			\begin{equation}\label{eq-CGSSF-(2)}
				\delta^{{\mathcal{H}}}\leq \frac{\left( r-2\right)}{2}\left\{ \frac{1}{\left(r-1\right) }\| \tau^{{\mathcal{H}}}\|^2 +\left( r+1\right) f_{1} \right\} +  \frac{(r-1)}{2} (3 f_2 - 2 f_3).
			\end{equation}
			The equality in $(\ref{eq-CGSSF-(2)})$ holds if and only if 
            $$\left\{\begin{array}{ll}
			\text{$r > 3$ and the statements \ref{S2}, \ref{S3}, \ref{S4}, \ref{S6}, and \ref{S7} follow}, & \text{when }~ f_3 < 0; \\
			\text{the statements \ref{S2}, \ref{S3}, \ref{S4}, and \ref{S7} follow}, & \text{when }~ f_3 = 0.
		\end{array} \right.$$
			
			\item Suppose $f_{2} \leq 0$ and $f_{3} > 0$, then
			\begin{equation}\label{eq-CGSSF-(3)}
				\delta^{{\mathcal{H}}}\leq \frac{\left( r-2\right)}{2}\left\{ \frac{1}{\left(r-1\right) }\| \tau^{{\mathcal{H}}}\|^2 +\left( r+1\right) f_{1} - 2 f_{3} \right\}.
			\end{equation}
			The equality in $(\ref{eq-CGSSF-(3)})$ holds if and only if the statements
            $$\left\{\begin{array}{ll}
			\text{\ref{S3}, \ref{S4}, \ref{S5}, and \ref{S8} follow}, & \text{when }~ f_2 < 0; \\
			\text{\ref{S3}, \ref{S4}, and \ref{S8} follow}, & \text{when }~ f_2 = 0.
		\end{array} \right.$$
			
			\item Suppose $f_{2}>0$ and $f_{3} > 0$, then
			\begin{equation}\label{eq-CGSSF-(4)}
				\delta^{{\mathcal{H}}}\leq \frac{\left( r-2\right)}{2}\left\{ \frac{1}{\left(r-1\right) }\| \tau^{{\mathcal{H}}}\|^2 +\left( r+1\right)  f_{1} +\frac{3(r-1)}{\left(r-2\right)}f_{2} -2f_{3}\right\}.
			\end{equation}
			The equality in $(\ref{eq-CGSSF-(4)})$ holds if and only if the statements $\ref{S2}$, $\ref{S3}$, $\ref{S4}$, $\ref{S7}$, and $\ref{S8}$ follow.
		\end{enumerate}
        In addition, it is understood that the term $f_3(r-1)$ in $(\ref{eq-CGSSF-(1)})$ and $(\ref{eq-CGSSF-(2)})$ vanishes when $f_3=0$.
	\end{theorem}
	
	\begin{proof}
		By (\ref{eq-GNSSF-(1)}), we have
		\begin{align}\label{eq-CGSSF-(5)}
			\delta^{\mathcal{H}} \leq & \left\{\frac{\left(r-2\right) }{2\left( r-1\right) }\|\tau^{\mathcal{H}}\|^{2} + \frac{1}{2}f_{1}\left( r^{2} - r - 2\right)  + \frac{3}{2} f_{2} \left( \|\mathcal{P}\|^{2}-2 \Theta(\mathbb{P})\right) + f_3 (1-r+\Phi(\mathbb{P}))\right\}.
		\end{align}
		We also simplify
		\begin{equation}\label{eq-CGSSF-(6)}
			\| \mathcal{P} \|^2 - 2 \Theta(\mathbb{P}) = \left\{\sum \limits_{j=3}^r (g_2(\pi_\ast h_1, \mathcal{P} \pi_\ast h_j))^2 + \sum \limits_{j=3}^r (g_2(\pi_\ast h_2, \mathcal{P} \pi_\ast h_j))^2+ \sum_{i,j=3}^{r} (g_{2}(\pi_{\ast} h_{i}, \mathcal{P} \pi_\ast h_{j}))^{2}\right\}.
		\end{equation}
		Clearly, from (\ref{eq-CGSSF-(6)}), we have $\| \mathcal{P} \|^2 - 2 \Theta(\mathbb{P}) \geq 0$. Also by the definitions, it follows that $\|\mathcal{P}\|^2 \geq 0$ and $\Theta(\mathbb{P}) \geq 0$. Further, we have $ 0 \leq \Phi(\mathbb{P}) \leq 1$ because
        \begin{equation}\label{eq_using_eta_xi}
            1 = \eta(\xi) = g_2(\xi, \xi) = \sum \limits_{j=1}^{r} g_2(\pi_\ast h_j, \xi) g_2(\pi_\ast h_j, \xi) = \sum \limits_{j=1}^{r} (\eta (\pi_\ast h_j))^2.
        \end{equation}
                We can also write $\mathcal{P} \pi_\ast h_j = \sum\limits_{i=1}^r g_2(\mathcal{P} \pi_\ast h_j, \pi_\ast h_i)\pi_\ast h_i$; therefore, $\|\mathcal{P} \pi_\ast h_j \|^2 = \sum \limits_{i=1}^r (g_2(\pi_\ast h_i, \mathcal{P} \pi_\ast h_j))^2$. Thus,
                \begin{equation}\label{norm_P_square2}
                    \sum \limits_{j=1}^{r} \|\mathcal{P} \pi_\ast h_j \|^2 = \| \mathcal{P} \|^2.
                \end{equation}
                Further, we also have
                \begin{align*}
                    \sum \limits_{j=1}^{r} \|\phi \pi_\ast h_j \|^2 & = \sum \limits_{j=1}^{r} g_2(\phi \pi_\ast h_j, \phi \pi_\ast h_j) \\&= \sum \limits_{j=1}^{r} g_2(\pi_\ast h_j, \pi_\ast h_j)- \sum \limits_{j=1}^{r}g_2(\pi_\ast h_j, \xi) g_2(\pi_\ast h_j, \xi) = r - g_2 (\xi, \xi)\\& = r-1 \quad \text{[since $g_2(\xi, \xi) = \eta(\xi) =1$]}
                \end{align*}
                and
                \begin{align*}
                    \| \phi \pi_\ast h_j \|^2 &= g_2(\phi \pi_\ast h_j, \phi \pi_\ast h_j) \\&= g_2(\mathcal{P} \pi_\ast h_j + \mathcal{Q} \pi_\ast h_j, \mathcal{P} \pi_\ast h_j + \mathcal{Q} \pi_\ast h_j) \quad [\text{using (\ref{decomposition_with_J})}] \\& = \| \mathcal{P} \pi_\ast h_j \|^2 + \| \mathcal{Q} \pi_\ast h_j \|^2.
                \end{align*}
                The above expressions mean $\sum \limits_{j=1}^{r}\| \mathcal{P} \pi_\ast h_j \|^2 + \sum \limits_{j=1}^{r} \| \mathcal{Q} \pi_\ast h_j \|^2 = r-1$; hence, (\ref{norm_P_square2}) gives
                \begin{equation*}
                    \| \mathcal{P}\|^2 \leq r-1.
                \end{equation*}
                Clearly, the maximum values of $\| \mathcal{P}\|^2$ and $\Phi(\mathbb{P})$ may reach $r-1$ and $1$, respectively. In addition, the minimum values of $\| \mathcal{P}\|^2$, $\Phi(\mathbb{P})$, and $\Theta(\mathbb{P})$ may attain $0$. Now, to maximize the right-hand side of (\ref{eq-CGSSF-(5)}), we have the following cases.
		\begin{enumerate}[$(1)$]
			\item Since $f_2 \leq 0$ and $f_3 \leq 0$, we must minimize $(\| \mathcal{P} \|^2 - 2 \Theta(\mathbb{P}))$ and $(1-r+\Phi(\mathbb{P}))$. Both attain their minimum values $0$ and $1-r$, respectively, when $(\| \mathcal{P} \|^2 - 2 \Theta(\mathbb{P})) = 0$ and $\Phi(\mathbb{P}) = 0$. In this case, (\ref{eq-CGSSF-(5)}) reduces to (\ref{eq-CGSSF-(1)}). In addition, if $f_2=0$ (resp. $f_3=0$), then $f_2(\| \mathcal{P} \|^2 - 2 \Theta(\mathbb{P}))$ (resp. $f_3(1-r+\Phi(\mathbb{P}))$) vanishes identically. Therefore, in the boundary cases, we do not need to minimize the respective terms, as they vanish.
			
			\item Since $f_{2}>0$ and $f_3 \leq 0$, we must maximize $\left(\| \mathcal{P} \|^{2} - 2 \Theta(\mathbb{P})\right)$ and minimize $\left(1 - r + \Phi (\mathbb{P} )\right)$, which attains the required values when $\| \mathcal{P} \|^2 =r-1$, $\Theta(\mathbb{P}) = 0$, and $\Phi(\mathbb{P})=0$. In this case, (\ref{eq-CGSSF-(5)}) reduces to (\ref{eq-CGSSF-(2)}). In addition, if $f_3=0$, then $f_3(1-r+\Phi(\mathbb{P}))$ vanishes identically; consequently the minimization of $\left(1 - r + \Phi (\mathbb{P} )\right)$ is not required.
			
			\item Since $f_2 \leq 0$ and $f_3 > 0$, we must minimize $(\| \mathcal{P} \|^2 - 2 \Theta (\mathbb{P}))$ and maximize $(1-r+\Phi(\mathbb{P}))$. Both achieve their required values $0$ and $(2-r)$, when $(\| \mathcal{P} \|^2 - 2 \Theta(\mathbb{P})) = 0$ and $\Phi(\mathbb{P}) = 1$. In this case, (\ref{eq-CGSSF-(5)}) reduces to (\ref{eq-CGSSF-(3)}). In addition, if $f_2=0$, then $f_2(\| \mathcal{P} \|^2 - 2 \Theta(\mathbb{P}))$ vanishes identically; consequently the minimization of $(\| \mathcal{P} \|^2 - 2 \Theta(\mathbb{P}))$ is not required.
			
			\item Since $f_{2}>0$ and $f_{3} > 0$, we must maximize the values of $\left(\| \mathcal{P} \|^{2}-2 \Theta (\mathbb{P} )\right)$ and $\left(1-r+\Phi (\mathbb{P} )\right)$, which are achieved when $\|\mathcal{P}\|^2=r-1$, $\Theta(\mathbb{P})=0$, and $\Phi(\mathbb{P})=1$. In this case, (\ref{eq-CGSSF-(5)}) reduces to (\ref{eq-CGSSF-(4)}).
		\end{enumerate}
		The equality conditions hold for the following reasons:
        \begin{itemize}
            \item The statements $\ref{S3}$ and $\ref{S4}$ are implied from the equality cases of $(\ref{eq-GNSSF-(1)})$ and $(\ref{eq-CGSSF-(5)})$, respectively.
            
            \item The condition $\| \mathcal{P} \|^2 - 2 \Theta(\mathbb{P}) = 0$ satisfies if and only if
            \begin{enumerate}[$(a)$]
                \item $g_2(\pi_\ast h_1, \mathcal{P} \pi_\ast h_j) =0$ and $g_2(\pi_\ast h_2, \mathcal{P} \pi_\ast h_j) =0$ for all $3 \leq j \leq r$, and
                \item $g_2 (\pi_\ast h_i, \phi \pi_\ast h_j) =0$ for all $3 \leq i, j \leq r$.
            \end{enumerate}
            The above conditions $(a)$ and $(b)$ are equivalent to the statement $\ref{S5}$.

            \item The condition $\Phi(\mathbb{P}) = 0$ satisfies if and only if $g_2(\pi_\ast h_1, \xi) = 0$ and $g_2(\pi_\ast h_2, \xi) = 0$. These are equivalent to the statement $\ref{S6}$.

            \item The condition $\Theta(\mathbb{P}) = 0$ satisfies if and only if $g_2(\pi_\ast h_1, \mathcal{P} \pi_\ast h_2) = 0$. This is equivalent to the statement $\ref{S2}$.

            \item The condition $\| \mathcal{P} \|^{2}=r-1$ satisfies if and only if $\mathcal{Q} \pi_\ast h_j = 0$ for each $1 \leq j \leq r$ due to the fact that $\sum \limits_{j=1}^{r}\| \mathcal{P} \pi_\ast h_j \|^2 + \sum \limits_{j=1}^{r} \| \mathcal{Q} \pi_\ast h_j \|^2 = r-1$. Then (\ref{decomposition_with_J}) yields
            \begin{equation*}
                \phi \pi_\ast h_j = \mathcal{P} \pi_\ast h_j \quad \text{for all $1 \leq j \leq r$}.
            \end{equation*}
                Equivalently, $\phi (\mathcal{R}) \subseteq (\mathcal{R})$; hence the leaves of $\mathcal{R}$ are invariant submanifolds of the target space \cite{HUAM}. Further, by \cite[Proposition 3.2]{Yano_Ishihara}, they are odd-dimensional almost contact submanifolds. Thus, the statement $\ref{S7}$ follows if and only if $\| \mathcal{P}\|^2 = r-1$.

            \item The condition $\Phi(\mathbb{P}) = 1$ satisfies if and only if $\left(\eta(\pi_\ast h_1)\right)^2 + \left(\eta(\pi_\ast h_2)\right)^2 = 1$. Using (\ref{eq_using_eta_xi}), we obtain that $\left(\eta(\pi_\ast h_1)\right)^2 + \left(\eta(\pi_\ast h_2)\right)^2 = 1$ if and only if
            \begin{equation*}
                g_2(\pi_\ast h_3, \xi) = g_2(\pi_\ast h_4, \xi) = \dots = g_2(\pi_\ast h_r, \xi) = 0.
            \end{equation*}
            The above expression is equivalent to the statement $\ref{S8}$.

            \item The conditions $\|\mathcal{P}\|^2 = r-1$, $\Theta(\mathbb{P})=0$, and $\Phi(\mathbb{P})=0$ cannot hold simultaneously when $r=3$. In fact, $\|\mathcal{P}\|^2 = r-1$ and $\Phi(\mathbb{P})=0$ force $\Theta(\mathbb{P})=1$ when $r=3$. Therefore, to satisfy the statements $\ref{S2}$, $\ref{S6}$, and $\ref{S7}$ simultaneously, we exclude $r = 3$.
        \end{itemize}
	\end{proof}
	
\begin{theorem}\label{thm_gssf_1_xi_perp}
		Let $\pi:(M^{m},g_{1})\to (N^{n}(f_{1},f_{2},f_{3}),g_{2}, \phi, \xi, \eta)$ be a Riemannian map from a Riemannian manifold to a generalized Sasakian space form with Reeb vector field $\xi \in \Gamma(\mathcal{R}^\perp)$. 
		\begin{enumerate}[$(1)$]
			\item Suppose $f_{2}\leq 0$, then
			\begin{equation}\label{eq-CGSSF-(1)_xi_perp}
				\delta^{{\mathcal{H}}}\leq \frac{\left( r-2\right)}{2}\left\{ \frac{1 }{\left(r-1\right) }\| \tau^{{\mathcal{H}}}\|^2 +\left( r+1\right) f_{1} \right\}.
			\end{equation}
			The equality in $(\ref{eq-CGSSF-(1)_xi_perp})$ holds if and only if the statements
            $$\left\{\begin{array}{ll}
			\text{\ref{S3}, \ref{S4}, and \ref{S5} follow}, & \text{when }~ f_2 < 0; \\
			\text{\ref{S3} and \ref{S4} follow}, & \text{when }~ f_2 = 0.
		\end{array} \right.$$
			
			\item Suppose $f_{2}> 0$, then
			\begin{equation}\label{eq-CGSSF-(2)_xi_perp}
				\delta^{{\mathcal{H}}}\leq \frac{\left( r-2\right)}{2}\left\{ \frac{1}{\left(r-1\right) }\| \tau^{{\mathcal{H}}}\|^2 +\left( r+1\right) f_{1} \right\} + \frac{3}{2} r f_2.
			\end{equation}
			The equality in $(\ref{eq-CGSSF-(2)_xi_perp})$ holds if and only if the statements \ref{S1}, \ref{S2}, \ref{S3}, and \ref{S4} follow.
		\end{enumerate}
	\end{theorem}
	
	\begin{proof}
		By (\ref{CFI_when_xi_in_R_perp}), we have
		\begin{equation}\label{eq-CGSSF-(5)_xi_perp}
			\delta^{\mathcal{H}} \leq \left\{\frac{\left(r-2\right) }{2\left( r-1\right) }\|\tau^{\mathcal{H}}\|^{2} + \frac{1}{2}f_{1}\left( r^{2} - r - 2\right) + \frac{3}{2} f_{2} \left( \|\mathcal{P}\|^{2}-2 \Theta(\mathbb{P})\right)\right\}.
		\end{equation}
        Using the fact $g_2 (\phi \pi_\ast h_j, \phi \pi_\ast h_j) = g_2(\pi_\ast h_j, \pi_\ast h_j)$ due to $\xi \in \Gamma(\mathcal{R}^\perp)$, one can easily show that
        \begin{equation*}
            \| \mathcal{P}\|^2 \leq r.
        \end{equation*}
        Then, similar to the previous theorem, we have the following cases.
		\begin{enumerate}[$(1)$]
			\item Since $f_2 \leq 0$, we must minimize $(\| \mathcal{P} \|^2 - 2 \Theta(\mathbb{P}))$, and it attains the minimum value $0$ when the right hand side terms in (\ref{eq-CGSSF-(6)}) vanishes. In this case, (\ref{eq-CGSSF-(5)_xi_perp}) reduces to (\ref{eq-CGSSF-(1)_xi_perp}). In addition, if $f_2=0$, then $f_2(\| \mathcal{P} \|^2 - 2 \Theta(\mathbb{P}))$ will be vanished identically; consequently the minimization of $(\| \mathcal{P} \|^2 - 2 \Theta(\mathbb{P}))$ is not required.
			
			\item Since $f_{2}>0$, we must maximize $\left(\| \mathcal{P} \|^{2} - 2 \Theta(\mathbb{P})\right)$, and it attains the maximum value when $\| \mathcal{P} \|^2 =r$ and $\Theta(\mathbb{P}) = 0$. In this case, (\ref{eq-CGSSF-(5)_xi_perp}) reduces to (\ref{eq-CGSSF-(2)_xi_perp}).
		\end{enumerate}
		In addition, the equality cases follow due to the equality case of (\ref{CFI_when_xi_in_R_perp}), (\ref{eq-CGSSF-(5)_xi_perp}), \cite[Proposition 3.1]{Yano_Ishihara}, and similar reasons mentioned in the previous theorems.
	\end{proof}
	
	Using the fact that $\tau^\mathcal{H}=0$ for harmonic Riemannian maps, we have the following corollary:
	\begin{corollary}\label{cor_harmonic_gssf}
		In particular, suppose $\pi:(M^{m},g_{1})\to (N^{n}(f_{1},f_{2},f_{3}),g_{2}, \phi, \xi, \eta)$ is a harmonic Riemannian map, and the Reeb vector field
		\begin{enumerate}[$(1)$]
			\item $\xi \in \Gamma(\mathcal{R})$, then
			\begin{equation*}
				\delta^{{\mathcal{H}}}\leq\left\{\begin{array}{ll}
					\frac{\left( r+1\right) \left(r-2\right)}{2} f_{1} - \left( r-1\right) f_3, & \text{if }~ f_{2}\leq 0 ~ \text{and} ~f_{3} \leq 0;\\
					\frac{\left( r+1\right) \left(r-2\right)}{2} f_{1} + \frac{3 r}{2} f_2 - \left( r-1\right) f_3, & \text{if }~ f_{2}> 0~ \text{and}~f_{3} \leq 0;\\
					(r-2)\left\{ \frac{\left( r+1\right)}{2}f_{1}-f_{3}\right\}, & \text{if }~ f_{2} \leq 0~ \text{and}~f_{3} > 0;\\
					(r-2)\left\{ \frac{\left( r+1\right)}{2}f_{1}-f_{3}\right\} + \frac{3r}{2}f_{2}, & \text{if }~ f_{2}>0 ~\text{and}~f_{3} > 0.
				\end{array} \right.
			\end{equation*}
			
			\item $\xi \in \Gamma(\mathcal{R}^\perp)$, then
			\begin{equation*}
				\delta^{{\mathcal{H}}} \leq 
				\left\{\begin{array}{ll}
					\frac{\left( r+1\right) \left(r-2\right)f_1}{2}, & \text{if }~ f_2 \leq 0;\\
					\frac{\left( r+1\right) \left(r-2\right) f_1 + 3r f_2}{2}, & \text{if }~ f_2 > 0.
				\end{array} \right.
			\end{equation*}
		\end{enumerate}
	\end{corollary}

 \section{A concluding remark}\label{sec_remarks}
     We note that the following listed results are motivated by analogs of certain results in submanifold theory; however, the proofs and expositions differ. 
     \begin{center}
         \begin{tabular}{cc}
        \hline
        \textbf{Result in current paper}  &  \textbf{Result(s) in submanifold theory}\\
        \hline
          Theorem $\ref{Theorem_GCSF_1}$ & $\cite[\text{\normalfont ~Theorem 3.3}]{KST}$\\ 
          
          Theorem $\ref{theorem_GSSF}$ & $\cite[\text{\normalfont ~Theorem 4.1}]{ACKY}$\\

          Item $(1)$ in Corollary $\ref{Corollary_GCSF}$ & $\cite[\text{\normalfont ~Lemma 3.2}]{Chen_1993}$\\

          Item $(2)$ in Corollary $\ref{Corollary_GCSF}$ & $\cite[\text{\normalfont ~Theorem 2}]{Chen_1996}$\\

          Item $(1)$ in Corollary $\ref{Corollary_1_GSSF}$ & $\cite[\text{\normalfont ~Theorem 3.2}]{Kim_Kim}$\\

          Item $(2)$ in Corollary $\ref{Corollary_1_GSSF}$ & $\cite[\text{\normalfont ~Theorem 3.2}]{GAH}$\\

          Item $(3)$ in Corollary $\ref{Corollary_1_GSSF}$ & $\cite[\text{\normalfont ~Theorem 3.2}]{Kim_Choe}$\\
          
          Theorem $\ref{theorem_f2geq0}$ &  $\cite[\text{\normalfont ~Theorem 3.2 and Theorem 3.3}]{ACKY}$\\
          
          Theorem $\ref{thm_gssf_1}$ & $\cite[\text{\normalfont ~Corollary 4.1}]{ACKY}$\\ 
          \hline
     \end{tabular}
     \end{center}
     In precise terms, for the expositions, our technique uses the general frameworks of Chen's inequalities and Riemannian maps. Thus, the established general Chen's first inequality provides an effective solution to establish Chen's first inequality, and hence $\delta$-invariants for Riemannian maps.

     \section*{Ethics approval}
     The work is original and is not submitted to other journals simultaneously.

     \section*{Funding}
     The author, R. Singh, is funded by the Human Resource Development Group (HRDG), Council of Scientific and Industrial Research (CSIR), New Delhi, India, through a research fellowship (File No: 09/0013(16054)/2022-EMR-I).

     \section*{Declaration of competing interest}
     The authors declare that they have no known competing financial interests or personal relationships that could have appeared to influence the work reported in this paper.

     \section*{Acknowledgements}
     We sincerely thank the anonymous referees for their constructive comments.


\begin{thebibliography}{99}
		\bibitem{ABC} Alegre, P., Blair, D.E., Carriazo, A.: Generalized Sasakian-space-forms. Israel J. Math., \textbf{141} (2004), 157-183.
		
		\bibitem{ACKY} Alegre, P., Carriazo, A., Kim, Y.H., Yoon, D.W.: B.-Y. Chen's inequality for submanifolds of generalized space forms. Indian J. Pure Appl. Math., \textbf{38} (2007), 185-201.
		
		\bibitem{ACM} Alegre, P., Chen, B.Y., Munteanu, M.I.: {Riemannian submersions, $\delta$-invariants, and optimal inequality}. Ann. Glob. Anal. Geom., \textbf{42}(3) (2012), 317-331.
		
		\bibitem{AEMMO} Arslan, K., Ezentas, R., Mihai, I., Murathan, C., \"Ozg\"ur, C.: Certain inequalities for submanifolds in $(\kappa, \mu)$-contact space forms. Bull. Austral. Math. Soc., \textbf{64} (2001), 201-212.

        \bibitem{Chen_totally_umbilical} Chen, B.Y.: {Totally umbilical submanifolds}. Soochow J. Math., \textbf{5} (1979), 9-37.
		
		\bibitem{Chen_1993} Chen, B.Y.: {Some pinching and classification theorems for minimal submanifolds}. Arch. Math., \textbf{60} (1993), 568-578.
		
		\bibitem{Chen_RIMA} Chen, B.Y.: {A Riemannian invariant and its applications to submanifold theory}. Results Math., \textbf{27} (1995), 17-26.
        
        \bibitem{Chen_1996} Chen, B.Y.: {A general inequality for submanifolds in complex-space-forms and its applications}. Arch. Math., \textbf{67} (1996), 519-528.
		
		
		
		\bibitem{Chen_2011} Chen, B.Y.: {Pseudo-Riemannian Geometry, $\delta$-invariants and Applications}. World Scientific, (2011).

        
		
		\bibitem{CDVV_1996} Chen, B.Y., Dillen, F., Verstraelen, L., Vrancken, L.: {An exotic totally real minimal immersion of $S^3$ in $\mathbb{C}P^3$ and its characterisation}. Proc. Royal Soc. Edinburgh, \textbf{126} (1996), 153-165.

        \bibitem{Chen_Verstraelen} Chen, B.Y., Verstraelen, L.: {A survey on quasi-umbilical submanifolds}. Paleis der Academi\"en, Brussel, (1976).
		
		\bibitem{Chen_2023} Chen, B.Y., V\^ilcu, G.E.: {Recent developments on the first Chen inequality in differential geometry}. Math., \textbf{11}(19) (2023), 1-50.
		
		\bibitem{Fischer_1992} Fischer, A.E.: {Riemannian maps between Riemannian manifolds}. Contemp. Math., \textbf{132} (1992), 331-366.

        \bibitem{GAH} Gupta, R.S., Ahmad, I., Haider, S.M.K.: B. Y. Chen's inequality and its applications to slant immersions into Kenmotsu manifolds. Kyungpook Math. J., \textbf{44} (2004), 101-110.

        \bibitem{HUAM} Hui, S.K., Uddin, S., Alkhaldi, A.H., Mandal, P.: Invariant submanifolds of generalized Sasakian-space-forms. Int. J. Geom. Methods Mod. Phys., \textbf{15}(9) (2018), 1-21.
		
		\bibitem{Kim_Choe} Kim, J.-S., Choi, J.: A basic inequality for submanifolds in a cosymplectic space form. Int. J. Math. Math. Sci., \textbf{9} (2003), 539-547.
		
		\bibitem{Kim_Kim} Kim, Y.H., Kim, D.-S.: A basic inequality for submanifolds in Sasakian space forms. Houst. J. Math., \textbf{25} (1999), 247-257.
		
		\bibitem{KST} Kim, J.-S., Song, Y.-M., Tripathi, M.M.: B.-Y. Chen inequalities for submanifolds in generalized complex space forms. Bull. Korean Math. Soc., \textbf{40}(3) (2003), 411-423.
		
		\bibitem{Meena} Meena, K., \c{S}ahin, B., Shah, H.M.: Chen's first inequality for Riemannian maps to complex space forms and $\delta$-invariants. Period. Math. Hungar., \textbf{91}(2) (2025), 588-603.

        \bibitem{Meulemeester_Verstraelen} Meulemeester, W.D., Verstraelen, L.: Codimension $2$ submanifold with a $2$-quasi-umbilical normal direction. J. Korean Math. Soc., \textbf{15}(2) (1979), 101-108.
		
		\bibitem{Nore_1986} Nore, T.: {Second fundamental form of a map}. Ann. Mat. Pura Appl., \textbf{146} (1986), 281-310.
		
		\bibitem{Olszak} Olszak, Z.: On the existence of generalized complex space forms. Israel J. Math., \textbf{65}(2) (1989), 214-218.
		
		\bibitem{Ozgur_2011} \"Ozg\"ur, C.: B.Y. Chen inequalities for submanifolds of a Riemannian manifold of quasi-constant curvature. Turkish J. Math., \textbf{35} (2011), 501-509.
		
		\bibitem{Ozgur_De} \"Ozg\"ur, C., De, A.: Chen inequalities for submanifolds of a Riemannian manifold of nearly quasi-constant curvature. Publ. Math. Debrecen, \textbf{82} (2013), 439-450.

       \bibitem{Roth_Upadhyay} Roth, J., Upadhyay, A.: Biharmonic submanifolds of generalized space forms. Differential Geom. Appl., \textbf{50} (2017), 88-104.
		
		\bibitem{Sahin_2010} \c{S}ahin, B.: {Invariant and anti-invariant Riemannian maps to K\"ahler manifolds}. Int. J. Geom. Methods Mod. Phys., \textbf{7}(3) (2010), 337-355.
		
		\bibitem{Sahin_2016} \c{S}ahin, B.: {Chen's first inequality for Riemannian maps}. Ann. Polon. Math., \textbf{117}(3) (2016), 249-258.
		
		\bibitem{Sahin_book} \c{S}ahin, B.: {Riemannian Submersions, Riemannian Maps in Hermitian Geometry, and Their Applications}. Elsevier, Academic Press, (2017).
		
		\bibitem{Tricerri} Tricerri, F., Vanhecke, L.: Curvature tensors on almost Hermitian manifolds. Trans. Amer. Math. Soc., \textbf{267}(2) (1981), 365-398.
		
		\bibitem{MMT_2003} Tripathi, M.M.: Certain basic inequalities for submanifolds in $(\kappa, \mu)$-spaces. Contemp. Math., \textbf{337} (2003), 187-202.

       \bibitem{MMT_conf_paper} Tripathi, M.M.: Different kinds of submanifolds of almost Hermitian manifolds. in: Proceedings of the 23rd International Differential Geometry Workshop on Submanifolds in Homogeneous Spaces \& Related Topics, (2021), 59-90.
		
		\bibitem{Vanhecke} Vanhecke, L.: Almost Hermitian manifolds with $J$-invariant Riemannian curvature tensor. Rend. Semin. Mat. Univ. Politec. Torino, \textbf{34} (1975-76), 487-498.
		
		\bibitem{Yano_Ishihara} Yano, K., Ishihara, S.: Invariant submanifolds of an almost contact manifold. Kodai Math. Sem. Rep., \textbf{21} (1969), 350-364.
        
        \bibitem{Yoon_2004} Yoon, D.W.: A basic inequality of submanifolds in quaternionic space forms. Balkan J. Geom. Appl., \textbf{9} (2004), 92-103.
	\end{thebibliography}
\end{document}